\documentclass{article}

\usepackage[letterpaper,top=2cm,bottom=2cm,left=3cm,right=3cm,marginparwidth=1.75cm]{geometry}

\usepackage{amssymb,amsfonts}
\usepackage{amsmath}
\usepackage{graphicx}
\usepackage{cite}
\usepackage[colorlinks=true, allcolors=blue]{hyperref}
\usepackage{mathtools}
\usepackage{xspace}
\usepackage{enumerate}
\usepackage{float}
\usepackage{xcolor}
\usepackage{multicol}
\usepackage{multirow}
\usepackage{pifont}
\usepackage{txfonts}
\usepackage{siunitx}
\usepackage[acronym,nomain]{glossaries}
\usepackage{url}
\usepackage{circuitikz}
\usepackage{orcidlink}
\usepackage{ragged2e}

\title{An Open Optimal Power Flow Model for the Australian National Electricity Market}

\begin{document}
\date{\vspace{-5ex}}
\maketitle

\centering 

\author{\centering Rahmat Heidari\footnote{R. Heidari and M. Amos are with the Energy Business Unit, Commonwealth Scientific and Industrial Research Organisation (CSIRO), Newcastle, NSW 2300, Australia, (e-mail: rahmat.heidarihaei@csiro.au, matt.amos@csiro.au)} \orcidlink{0000-0003-4835-2614},{} 
        Matthew Amos\footnotemark[1], 
        Frederik Geth\footnote{F. Geth is with the GridQube, Brisbane, QLD 4300, Australia, (e-mail: frederik.geth@gridqube.com)} \orcidlink{0000-0002-9345-2959}{}
}

%

\raggedleft
\justifying


\begin{abstract}
    The Australian National Electricity Market (NEM) is a complex energy market that faces challenges due to the increasing number of distributed energy resources (DERs) and the transition to a net-zero emissions target. Power system modelling plays a crucial role in addressing these challenges by providing insights into different scenarios and informing decision-making. However, accessing power system data containing sensitive information can be a concern. Synthetic data offer a solution by allowing researchers to analyze and develop new methods while protecting confidential information. This paper utilizes an existing synthetic network model based on the NEM (`S-NEM2300'-bus system) to develop a benchmark for power system optimization studies. The model is derived and enhanced using PowerModels.jl and MATPOWER data models, and feasibility is ensured through power flow and optimal power flow studies. The resulting benchmark model, called `S-NEM2000'-bus system, is validated and enriched with additional parameters such a thermal limits, generation fuel categories and cost models. The `S-NEM2000'-bus system is an \emph{open}\footnote{\url{https://github.com/csiro-energy-systems/Synthetic-NEM-2000bus-Data}} dataset which provides a valuable resource for optimization studies in the power system domain. 

\end{abstract}

\paragraph{Disclaimer:}
The open network data developed in this study, specifically the `S-NEM2000'-bus benchmark model, is a representation created for research purposes and does not reflect the actual NEM. It is important to note that this model is a work in progress and may not include all the features and components present in the real NEM. Furthermore, the model can be further improved by incorporating additional features and elements and refining existing components. The ongoing development and growth of this work aim to enhance the representation of the NEM and its associated power system for future studies and analyses. The authors  invite collaborators to reach out through the issue tracker in the \href{https://github.com/csiro-energy-systems/Synthetic-NEM-2000bus-Data}{repository}  for any inquiries, feature requests, or collaboration opportunities.


\section{Introduction}
\label{sec:introduction}

		

\begin{figure}[htbp]
  \centering
  \includegraphics[width=0.5\columnwidth]{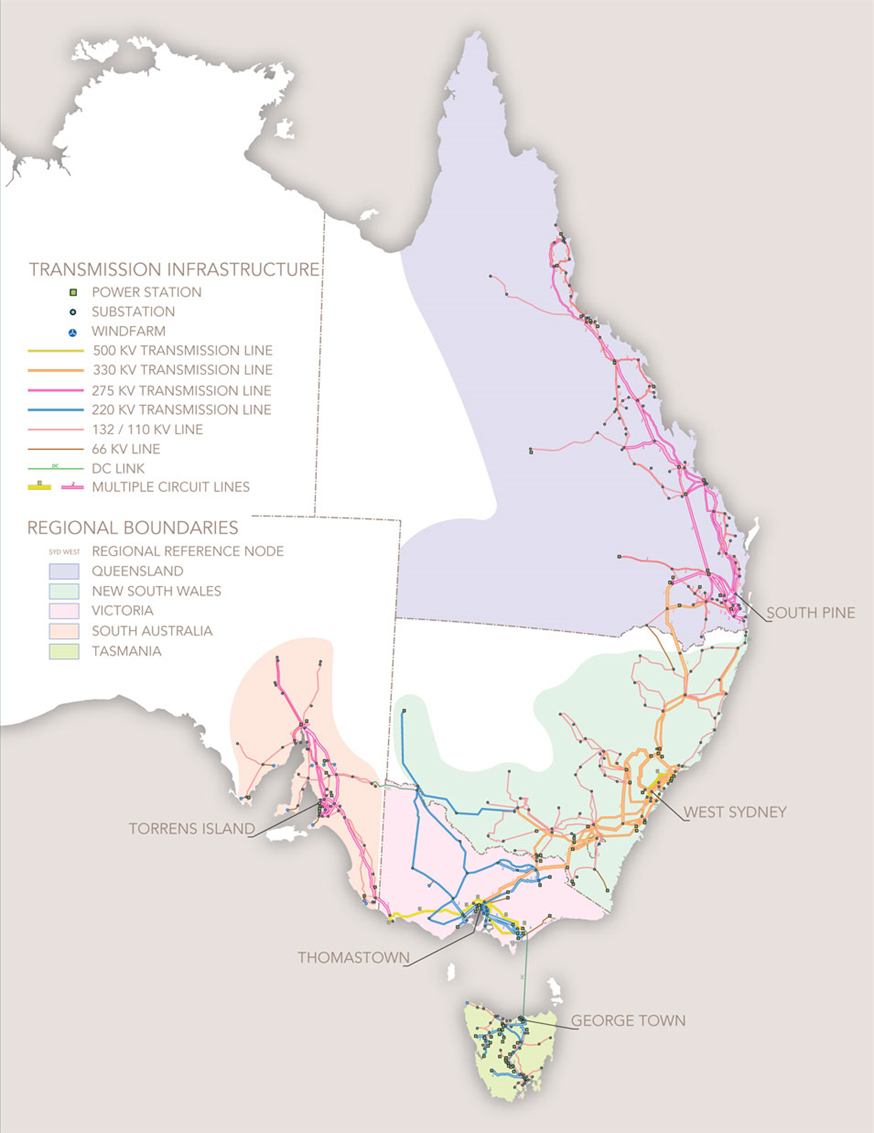}
  \caption{One-line diagram of the Australian NEM  \cite{AEMC2018}.}
  \label{fig:NEM}
\end{figure}

The Australian National Electricity Market (NEM) is a complex energy market that operates across five states and territories in Australia and regulates and manages the supply and demand of electricity within the market. The NEM, shown in Figure~\ref{fig:NEM}, is one of the largest interconnected power systems in the world, covering 40,000 km of transmission lines and cables, and servicing around 10 million customers\footnote{\url{https://www.aemo.com.au/-/media/Files/Electricity/NEM/National-Electricity-Market-Fact-Sheet.pdf}}. The market operates on a real-time basis, balancing the demand for electricity with the supply of electricity from a range of sources, including large-scale generators, renewable energy systems, and distributed energy resources (DERs).

The NEM is facing challenges due to the rising number of DER such as rooftop solar, batteries, and electric vehicles. Challenges include integrating DERs into the existing grid infrastructure, ensuring system stability and reliability, and addressing grid congestion in areas with high penetration of rooftop solar, which can cause voltage fluctuations and result in power outages \cite{Arrano_Vargas15020429}.  This requires significant investment in network infrastructure, as well as changes to the regulatory framework to support the integration of DERs.
The transition to a net-zero emissions target in the NEM also presents challenges, including the need for energy market reform and changes to regulatory frameworks to incentivize renewable energy resources deployment and emissions reduction. Additionally, the uncertainty surrounding government policy and support for renewable energy can make it difficult for investors to make long-term investment decisions.

Power system modelling is essential in dealing with the above-mentioned challenges in the electricity sector. Modelling tools provide insights into the impacts of different scenarios, such as the deployment of renewable energy sources or grid upgrades, and inform decision-making related to planning, design, operation, policy, regulation, and investment. This helps to ensure the reliability, sustainability, and affordability of the power system.

Despite necessity of access to power system data that enables reliable modelling and decision-making,   power system data can contain confidential information such as operational parameters, network topology, and customer data.   By using \emph{synthetic data}, this sensitive information can be protected while still allowing researchers and engineers to analyse and develop new methods.
Synthetic data are beneficial in the power systems context as they allow engineers to simulate and analyse power systems under different scenarios, identify potential issues, and test new algorithms and techniques. This helps to ensure the reliability and safety of the power system, which is critical to the functioning of society as a whole.

Furthermore, openly accessible data, \emph{open data}, plays an important role in power systems as it contributes to increased transparency, efficiency, and innovation. Open data encourages collaboration and knowledge sharing among different organizations and researchers, leading to the development of best practices and the advancement of the power sector as a whole.

%
%

	The authors of \cite{Coffrin2014, PGlib2021} provided a comprehensive library of all publicly available transmission network datasets (to date of the paper) under a common data format to be used for AC power flow optimisation studies.  To this end, the authors updated the dataset by reconstructing the missing data and adding key parameters such as generator cost functions and thermal limits, for algorithmic benchmarking purposes. 
    More recently,  a detailed network model of the European countries has been released in \cite{TOSATTO2022112907} to study market simulations of different combinations of wind and transmission capacity that is installed in the North Sea. The modelling framework and datasets are open which facilitates a range of other studies to investigate the impact of large-scale projects on the European system and electricity markets. 
    An open grid data of the NEM transmission network has been developed and released in \cite{Xenophon2018} which includes network and generation data sets, geospatial locations of network elements, and is designed to interface with AEMO's public database. This interface is specifically advantageous since it allows historic regional load and generator dispatch data to be integrated with the network model.  The authors ensured completeness and functionality of the dataset by assessing power flow models.

    Recently a synthetic network model based on the NEM has been released \cite{Arrano-Vargas2022, ArranoVargas2022, Arrano-Vargas2021}. This model, called S-NEM2300-bus benchmark model,  has longitudinal structure and is useful for EMT real-time digital simulations.  The authors developed this model based on a PSS/E model of the NEM and used statistical methods to obfuscate generation, transmission and load parameters.  The resulting model is a relational database released as CSV files and can be used in real-time EMT platforms.

    In this work we use this S-NEM2300 network data to develop a benchmark for power system \emph{optimisation} studies. We develop a data model to represent ac-feasible steady state power flow behaviour of the synthetic NEM model by parsing, cleaning and enhancing the input relational data in csv format (extract from Hypersim simulations) to PowerModels.jl \cite{8442948} and MATPOWER \cite{5491276} data models.
	PowerModels.jl \cite{8442948} is a Julia/JuMP package to study steady-state power network optimization problems.  It provides capabilities such as parsing and modifying network data as well as a common platform for computational evaluation of emerging power network formulations and algorithms.  PowerModels' network data model is organized as data dictionaries and the key names are designed  consistent with MATPOWER's \cite{5491276} file format which is a familiar data model for power system researchers, with some exceptions on the having separate components for loads and shunts, and design of time series simulation as ``multi-network''.

    To ensure feasibility of the converted data we conduct a series of power flow and optimal power flow studies as described in Section~\ref{sec:methodology}. The ac-feasible converted data is then validated by comparing power flow results with the steady state results in the original data release.  The converted synthetic network data is then complemented and extended by adding thermal limits and generation fuel type and cost models to serve as a benchmark model for optimisation studies.

\subsection{Motivation}
\label{sec:motivation}
The S-NEM2300 network is developed as a benchmark model based on the NEM which is compatible for real-time simulations. We aim to convert this data to a useful data for optimal power flow studies. The original data comes with the following limitations:
\begin{itemize}
    \item It is generated based on a snapshot of the actual NEM at a Summer day of 2018. As such, it is likely that not all generators were scheduled at the time, hence not represented in the synthetic model. 
    \item Generators are classified as thermal or hydro, but further details regarding fuel type and cost data are not taken into account. 
    \item Line and transformer thermal ratings are not identified.
\end{itemize}

The main motivation to develop a benchmark model is studying the following applications on the NEM-level network:
\begin{itemize}
    \item We aim to develop the data model to represent ac-feasible steady state power flow behaviour of the S-NEM2300 model. That is, the AC formulations are feasible for the developed network data, which ensures feasibility of other formulations such as current-voltage (IV), Quadratic Convex (QC), and `DC'.
    \item Impact of DER integration on the grid can be studied in a large scale network with similar properties the the NEM.
    \item Conducting optimal power flow studies such as unit commitment, economic dispatch, and market studies.
    \item Conducting security constrained optimal power flow studies.
    \item Development of hybrid ac-dc network by adding the NEM HVDC lines in the synthetic model and investigating security of the network against credible contingencies.
\end{itemize}

\section{Original Synthetic NEM (S-NEM2300) Network Data}
\label{sec:original_network_description}
\begin{figure}[H]
  \centering
  \includegraphics[width=0.4\columnwidth]{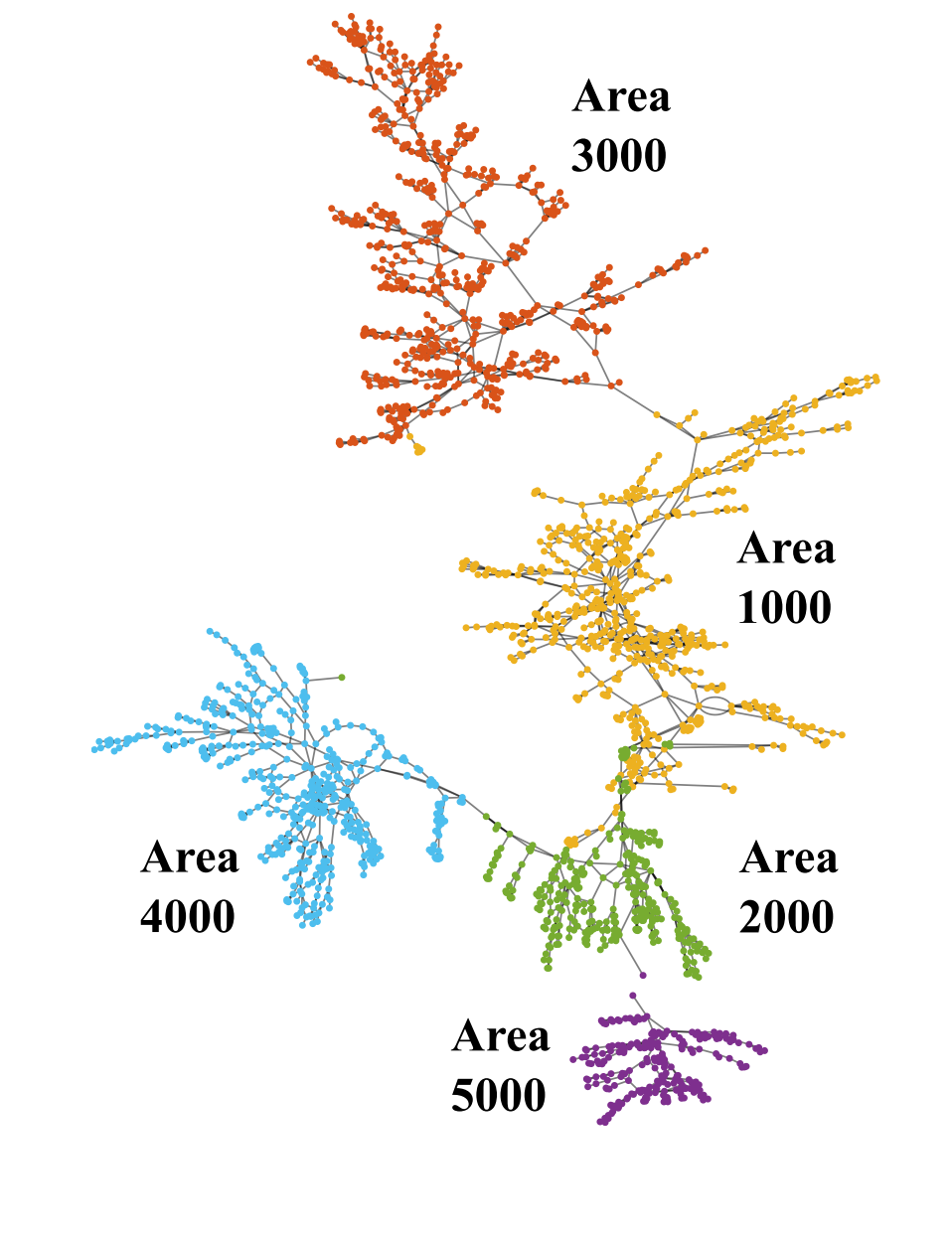}
  \caption{One-line graph diagram of the synthetic NEM data \cite{ArranoVargas2022}.}
  \label{fig:NEM_synthetic}
\end{figure}

We first revisit the main properties of the synthetic NEM data released in  \cite{ArranoVargas2022}. Table~\ref{table_original_nem_elements} and Figure~\ref{fig:NEM_synthetic} illustrate the main network components and the areas of the network based on synthetic coordinates available in the dataset. More properties of the network such as area interconnectors,  graph properties, structural properties and area diameters are provided in \cite{ArranoVargas2022}.

\begin{table}[tbh]
  \centering 
   \caption{Original synthetic NEM network number of components and areas \cite{ArranoVargas2022}.} 
   \label{table_original_nem_elements}
    \begin{tabular}{l  c c }
    \hline
    Component & Count \\
    \hline
    Buses & 2340 \\
    Lines & 1821 \\
    Two-winding transformers & 1418 \\
    Three-winding transformers & 113 \\
    Generators & 265\\
    Shunts & 301 \\
    Loads & 1696 \\
\hline
\end{tabular}
\hspace{12ex}
\begin{tabular}{l  c c }
\hline
    Area & State & Bus ID range \\
    \hline
    1000 & New South Wales & 1000-1999 \\
    2000 & Victoria & 2000-2999 \\
    3000 & Queensland & 3000-3999 \\
    4000 & South Australia & 4000-4999 \\
    5000 & Tasmania & 5000-5999 \\
    \\
    \\
\hline
\end{tabular}
\end{table}

%

\section{Data Model Development Framework}
\label{sec:data_model_development}
This section describes data model cleaning and derivation process and the methodology to develop an ac-feasible network data.  The input data is in CSV format which we first parse into PowerModels data model that can then readily be serialized into MATPOWER data format\footnote{command \texttt{export$\_$matpower(file, data)}}.

\subsection{Data Models}
\label{sec:data_models}

\subsubsection{Hypersim}
Opal-RT Hypersim\footnote{\url{https://www.opal-rt.com/systems-hypersim/}} is a powerful real-time simulation platform that can be used for power system simulations. It is a software solution that allows engineers to create a virtual environment that mimics the physical power system, enabling accurate modeling and testing of complex power systems.

\subsubsection{MATPOWER}
The MATPOWER data model format is a standardized data format used by the MATPOWER software package \cite{5491276} for representing power system data in a MATLAB file.  It provides an established, flexible format that researchers know, which makes it easier to analyse, compare, and share data in the power system research community.

\subsubsection{PowerModels.jl}
The PowerModels.jl package internally uses a dictionary to store network data with strings as keys so that it can be serialized to JSON for algorithmic exchange \cite{8442948}. A detailed description of the PowerModels data model is provided in the documentation\footnote{\url{https://lanl-ansi.github.io/PowerModels.jl/stable/}}.
PowerModels' data model can be instantiated by reading MATPOWER case studies, and PowerModels can serialize data to both MATPOWER and JSON.




\subsection{Data Model Derivation}
\label{sec:data_model_derivation}
We parse and clean the input relational data into the PowerModels data model with bus, branch, load, gen, switch, storage, shunt and dcline as dictionary keys.  Each network element is indexed and the data is stored in per unit values.  Transformers are stored as branches with an identifier flag. 

The data cleaning process consists of the following steps:
\begin{itemize}
	\item Identifying how to accurately convert each element data from the relational dataset to the PowerModels format.   For most network elements such as buses, lines, shunts, load and generators, the cleaning is fairly straightforward.  For two-winding and three-winding transformers, however, the cleaning requires implementing an equivalent circuit, which explained later in this section.  Note that the input data does not contain any dc line, so in this paper we follow the same approach, and develop a model with dc line and buses in a future work.
	
	\item Regional areas and interconnector branches are identified and flagged.
		
	\item Per unit conversion and unit transformations: We convert all data to per unit w.r.t their individual base values as given in the original data set (per unit conversion) and then update the per unit values to be on the base value of the cleaned data model (unit transformation).
	
\end{itemize}

\subsubsection{Two-Winding Transformer}
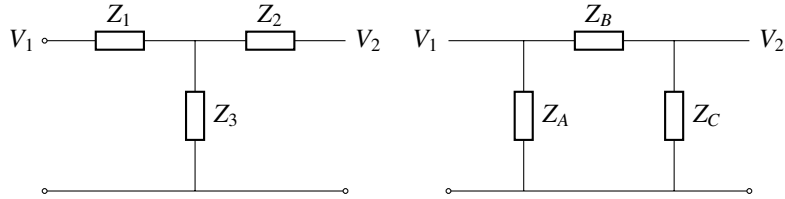
\begin{figure}[H]
  \centering
  
\ctikzset{bipoles/length=.8cm}
\begin{circuitikz}[european, scale=0.5] \draw 
  (0,0) to[short,o-o] (8,0)
  (0,4)  node[anchor=east] {$V_1$} to[R=$Z_1$,  o-] (4,4)  to[R=$Z_2$] (8,4)   node[anchor=west] {$V_2$}
  (4,4) to[R=$Z_3$] (4,0) ;
\end{circuitikz}
\ctikzset{bipoles/length=.8cm}
\begin{circuitikz}[european, scale=0.5] \draw 
  (0,0) to[short,o-o] (8,0)
  (0,4)  node[anchor=east] {$V_1$}  -- (2,4) to[R=$Z_B$] (6,4)  -- (8,4)   node[anchor=west] {$V_2$}
  (2,4) to[R=$Z_A$] (2,0) 
  (6,4) to[R=$Z_C$] (6,0) ;
\end{circuitikz}

  \caption{Two winding transformer equivalent circuits: T model (left) and $\Pi$ model (right).}
  \label{fig:2winding_tr_equivalent}
\end{figure}

In the input data, two-winding transformer models include primary, secondary, and magnitization impedance values. This model represents a T-section circuit as shown in Figure~\ref{fig:2winding_tr_equivalent} on the left. The branch model in PowerModels data format, however, is in $\Pi$-section Figure~\ref{fig:2winding_tr_equivalent} on the right. We therefore perform a T to $\Pi$ circuit conversion using equations in \eqref{eq:T-Pi-conversion} and calculate the transformer tap ratio from primary and secondary per unit voltage values.

\begin{align}
	\label{eq:T-Pi-conversion}
	Z_A = \frac{Z_1 Z_2 + Z_2 Z_3 + Z_1 Z_3}{Z_2}, \quad
	Z_B = \frac{Z_1 Z_2 + Z_2 Z_3 + Z_1 Z_3}{Z_3}, \quad
	Z_C = \frac{Z_1 Z_2 + Z_2 Z_3 + Z_1 Z_3}{Z_1}
\end{align}


\subsubsection{Three-Winding Transformer}
\begin{figure}[H]
  \centering
\ctikzset{bipoles/length=.8cm}
\begin{circuitikz}[scale=0.5] \draw 
  (0,0) to[short,o-o] (10,0)
  (0,4)  node[anchor=east] {$V_1$} to[R=$R_1$,  o-] (2,4) to[L=$L_1$] (4,4) -- 
  		(6,4) -- (6,5) to[R=$R_2$] (8,5) to[L=$L_2$, -o] (10,5)   node[anchor=west] {$V_2$}
  		(6,4) -- (6,3) to[R=$R_3$] (8,3) to[L=$L_3$, -o] (10,3)  node[anchor=west] {$V_3$}
  (4,4) -- (4,3) -- (3,3) to[R=$R_0$] (3,1) -- (4,1) -- (4,0)
                (4,3) -- (5,3) to[L=$L_0$] (5,1) -- (4,1)
;
\end{circuitikz}
  \caption{Three winding transformer equivalent circuit.}
  \label{fig:3winding_tr_equivalent}
\end{figure}
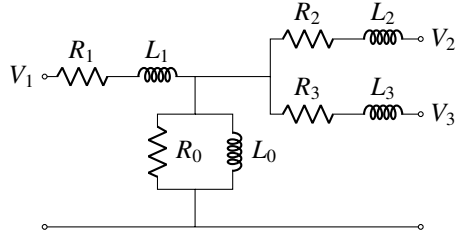

An equivalent circuit of a three-winding transformer is shown in Figure~\ref{fig:3winding_tr_equivalent}.  This circuit is implemented in PowerModels format as three branches representing three two-winding transformers that are connected to a common auxiliary node in between. The primary, secondary and tertiary impedance values are available from input data and the magnetization impedance is considered to be at the primary side,  i.e. the primary two winding transformer.
To accurately capture the actual three winding transformer tap ratios, each two-winding transformer tap ratio is calculated from its per unit voltage divided by the auxiliary node per unit voltage which is assumed to be 1. 


\subsection{Methodology}
\label{sec:methodology}
While the original network data is implemented and tested in Hypersim, cleaning and derivation of a different data model is still challenging as the cleaned data may cause infeasibility in optimal power flow (OPF) studies, as OPF needs to satisfy bounds such as voltage and current limits, and avoid voltage collapse.
A technique to analyse power flow infeasibility is to extend traditional power-balance constrained OPF to include nodal slack variables that indicate constraint violations when there is no feasible solution.  
The infeasibility analysis technique is such that if the network is feasible, the slack variables have values of zero at convergence. Otherwise, nonzero slack variables can be used to identify and localize infeasibility for planning insight or corrections during operation.
This method was first introduced in \cite{336130} for transmission systems and then extended in other work such as \cite{8960553, li2019lassoinspired}. Lately this technique has been applied to distribution systems in \cite{foster2023actionable}.

\begin{figure}[htbp]
  \centering
    \includegraphics[width=0.7\columnwidth]{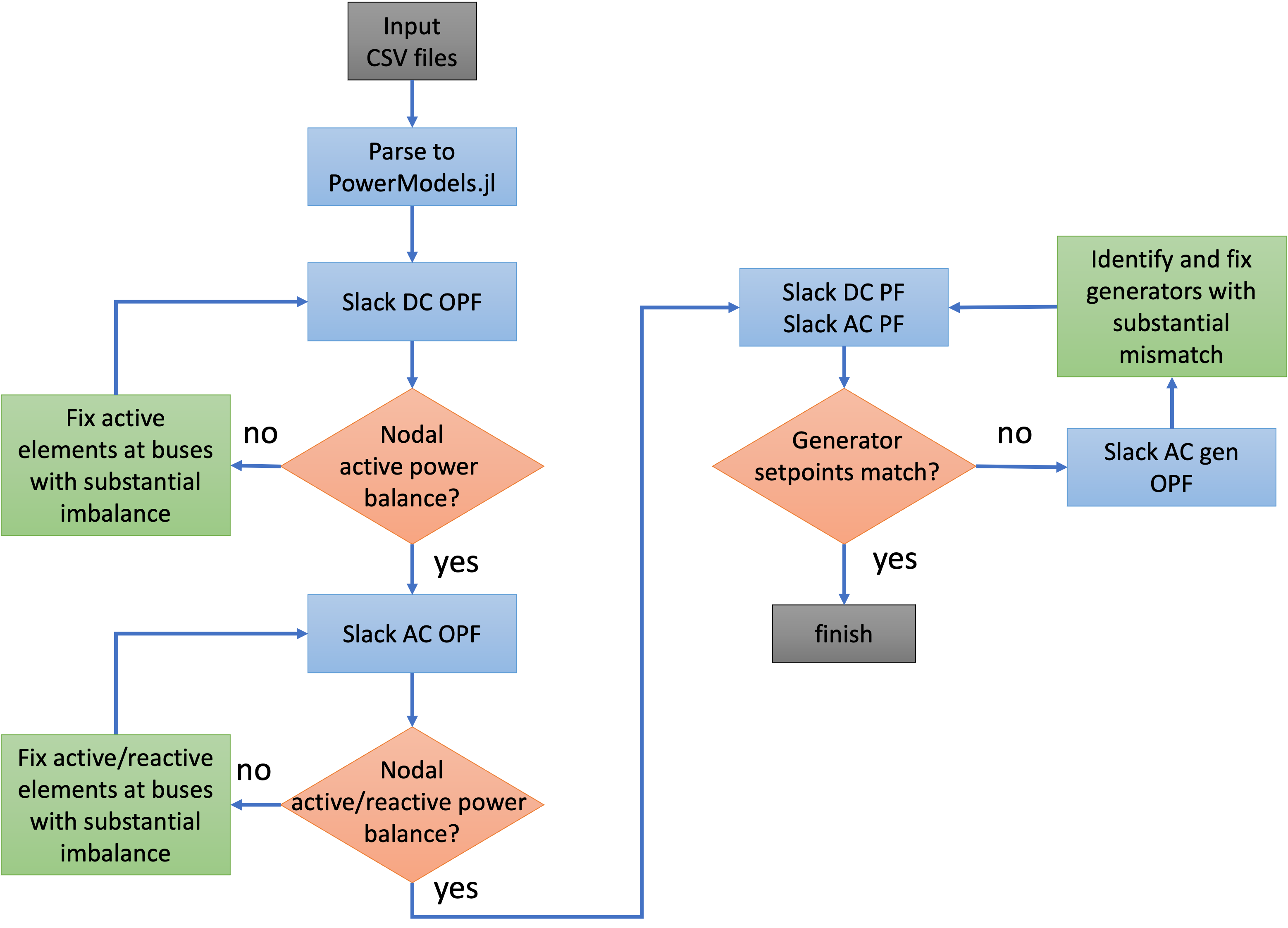}
  \caption{Model conversion algorithm consisting of a series of relaxed power flow and optimal power flow evaluations.}
  \label{fig:methodology}
\end{figure}

In this work we apply the infeasibility analysis method used in \cite{li2019lassoinspired} in an algorithmic way to identify different local constraint violations and the physical meaning behind it. Figure~\ref{fig:methodology} shows the flowchart of the algorithm implemented to identify infeasibility issues in the derivation process.  Starting from input files, the cleaned data is first tested by running a linearized `DC' OPF with slack nodal power variables to identify whether active power balance at all nodes are maintained. Nodes with significant imbalance are investigated as to what active elements could be the problem. The next step is to identify reactive power balance at the nodes which is the outcome of an AC-OPF with nodal slack power variables. Similar to the previous step, the nodes with substantial imbalance are investigated and the derivation process is modified until active and reactive power balance at all nodes are maintained. If the slack AC-OPF renders perfectly balanced network, the derivation could conclude at this point, but in case of existing imbalance in this step, the next two steps prove effective. The next step is to test generation dispatch setpoints, which is achieved by comparing AC and linearized `DC' PF solutions of the network with slack nodal power variables. If the generator setpoints mismatch, slack power variables are introduced to generators to identity which generators cause the most mismatch and the derivation process is modified accordingly.

\subsection{Validation of the Derived Synthetic NEM (S-NEM2000) Network Data}
\label{sec:validation}
\begin{figure}[htbp]
  \centering
    \includegraphics[width=1\columnwidth]{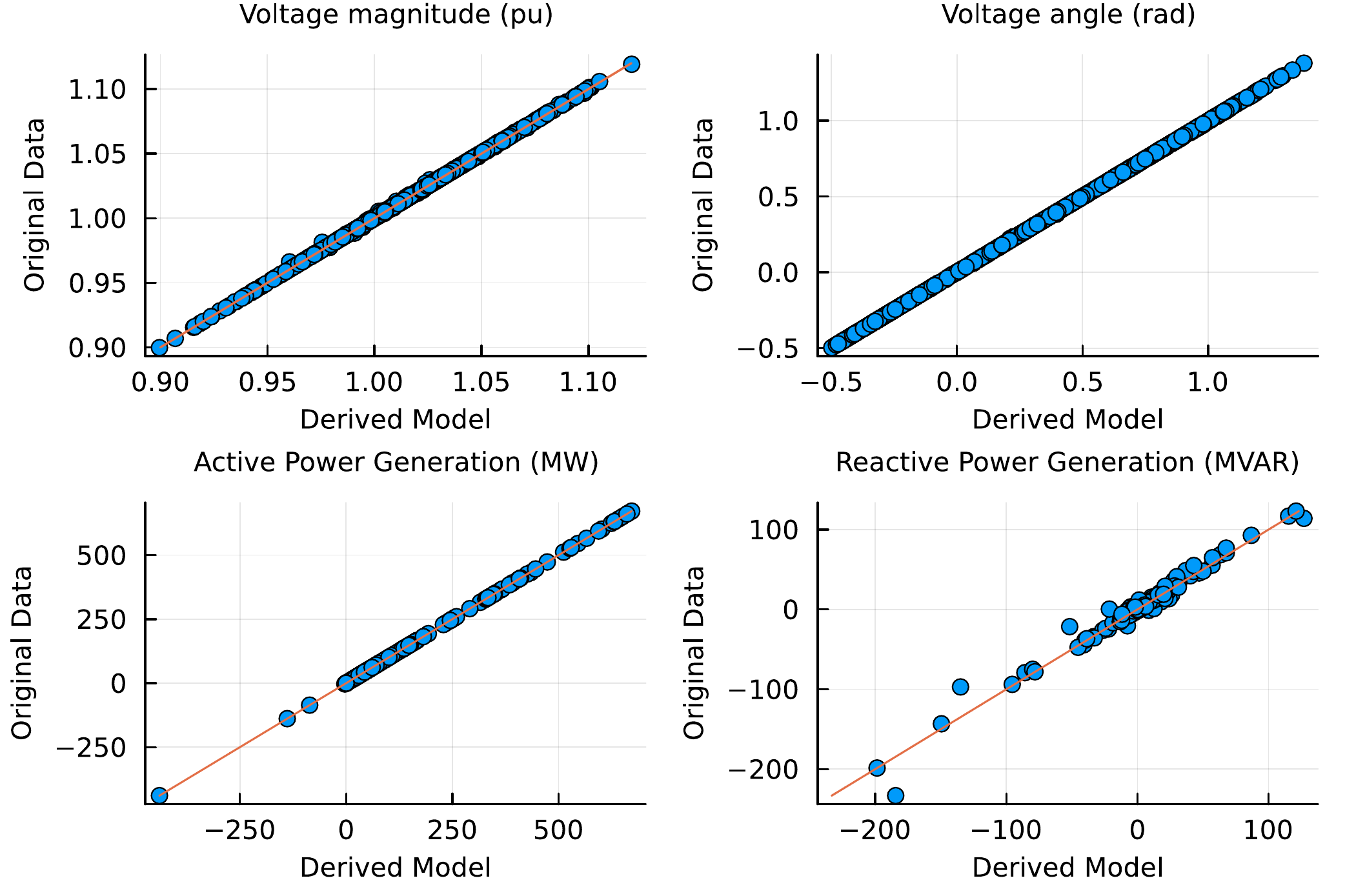}
  \caption{Validation of the derived model power flow results.}
  \label{fig:validation_base}
\end{figure}

Figure~\ref{fig:validation_base} demonstrates how the power flow results of the derived network data compares to the steady state results of the original network. As can be seen, the derived model accurately tracks the original model in voltage magnitude and phase angle, and generators' active power setpoints. The slight mismatch in generators' reactive power setpoints is mainly due to approximations made in generator models which  ignore internal impedance values.

\section{Thermal Limits}
\label{sec:thermal_limits}
The S-NEM2300 network data was developed to serve for phasor and EMT simulation studies. To enhance the network data to be useful for \emph{optimisation} studies we need to add reasonable network constraints such as line and transformer thermal limits.  

We adopt the data driven approach described in \cite{Coffrin2014} built on publicly available datasets. As noted in \cite{Coffrin2014},  these datasets reflect many of the statistical features that are common on real networks.  Section~\ref{sec:a_thermal_limit_model} provides an overview of the reference datasets and the derivation and validation of the thermal limit models.

To apply the thermal limit models to the synthetic NEM model, however,  we first review the  properties of the S-NEM2300 data.  The network data is then modified to reasonably resemble the transmission systems by reducing network through joining short lines, and removing ideal lines and circuit breakers.  The data driven thermal limit models are then applied to the modified synthetic NEM network data which we call S-NEM2000.  We next run an algorithm that links and updates thermal limits of the connecting lines and transformers so that it more closely represents a real power system set of constraints.  The properties of the modified network and linked thermal limits are validated and discussed.

\subsection{A Thermal Limit Model }
\label{sec:a_thermal_limit_model}
In practice,  line thermal limits are straightforward to calculate if conductor type and line length are available. In this synthetic dataset, however, we only have line impedance values and the connecting bus nominal voltages.  With this limited information we use a linear regression model proposed in \cite{Coffrin2014} to determine line thermal limits based on the values of line reactance $X$,  line resistance $R$ and nominal voltage $v$.  For the elements that any of these values are not available or not applicable (in case of ideal lines and transformers),  a reasonable upper bound is derived.

\subsubsection{Reviewing the Reference Data}
\label{sec:thermal_limits_ref_data}
\begin{figure}[H]
  \centering
  \includegraphics[width=0.48\columnwidth]{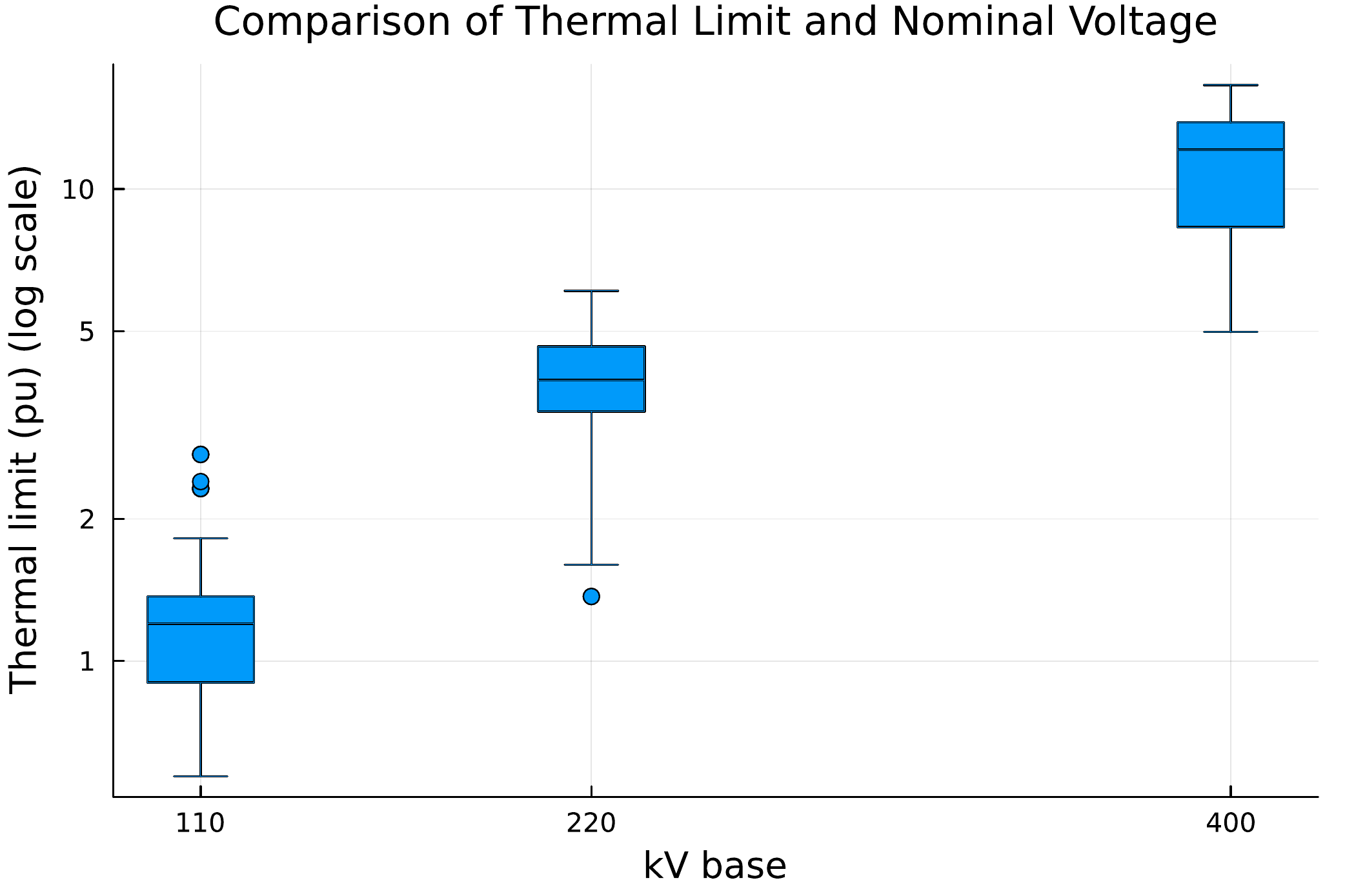}
    \includegraphics[width=0.48\columnwidth]{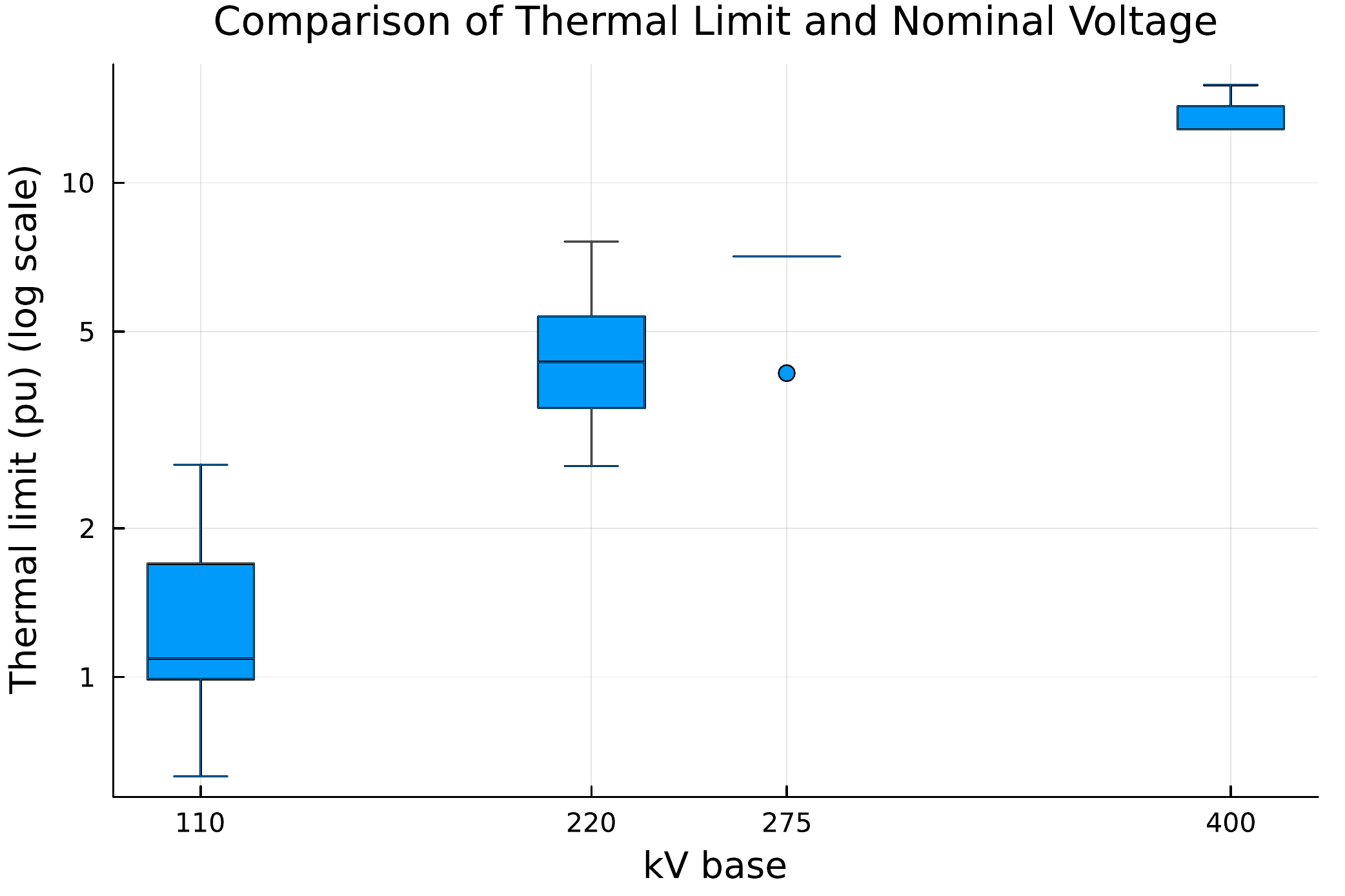}
  \caption{Polish network (left) and EIRGrid network (right) thermal limits and voltage levels correlation.}
  \label{fig:polish_irish_thermal_limit_kv_base}
\end{figure}

We use two sources of transmission thermal limits for this benchmark study: the Polish transmission system data \cite{5491276}, and the Irish transmission transmission system data \cite{EIRGrid2020}.   Although some transformer data is also available in these datasets,  we focus the data driven study only on the available transmission line data.  In order to stay within reasonable range of thermal limits and line types, we only consider lines with $X/R$ ratio below $20$ and normalised thermal rating (thermal rating divided by nominal voltage,  that is, current rating) above $0.005$.  Figure~\ref{fig:polish_irish_thermal_limit_kv_base} shows the correlation between line thermal limits and the nominal voltage.  As stated in \cite{Coffrin2014}, the significant variance at each voltage level which even overlaps with the neighbouring voltage levels replicates a realistic behaviour in real networks where a given nominal voltage may contain low and high capacity lines.

\subsubsection{Statistical Model}
\label{sec:statistical_model}
\begin{figure}[htbp]
  \centering
  \includegraphics[width=0.48\columnwidth]{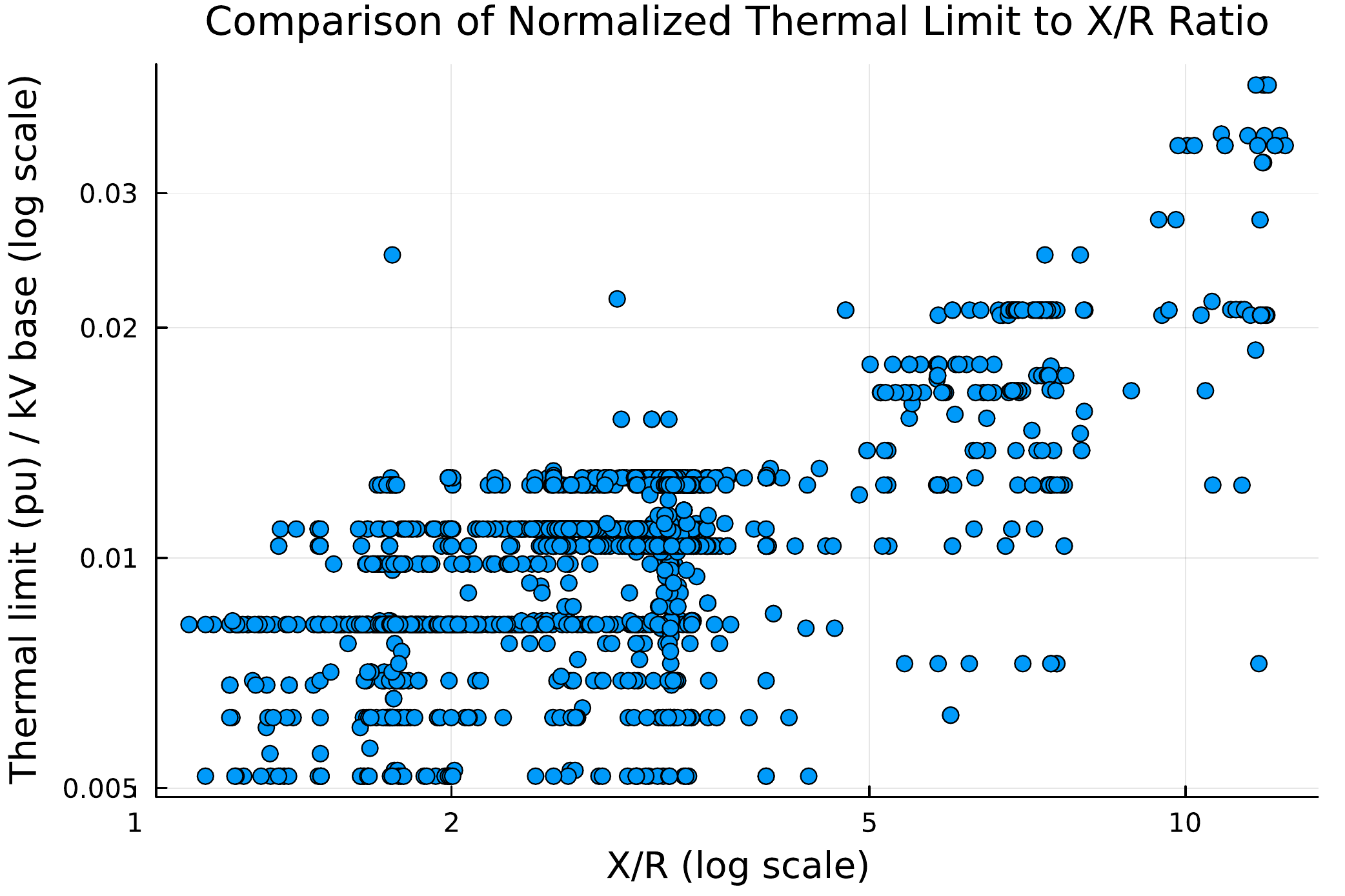}
    \includegraphics[width=0.48\columnwidth]{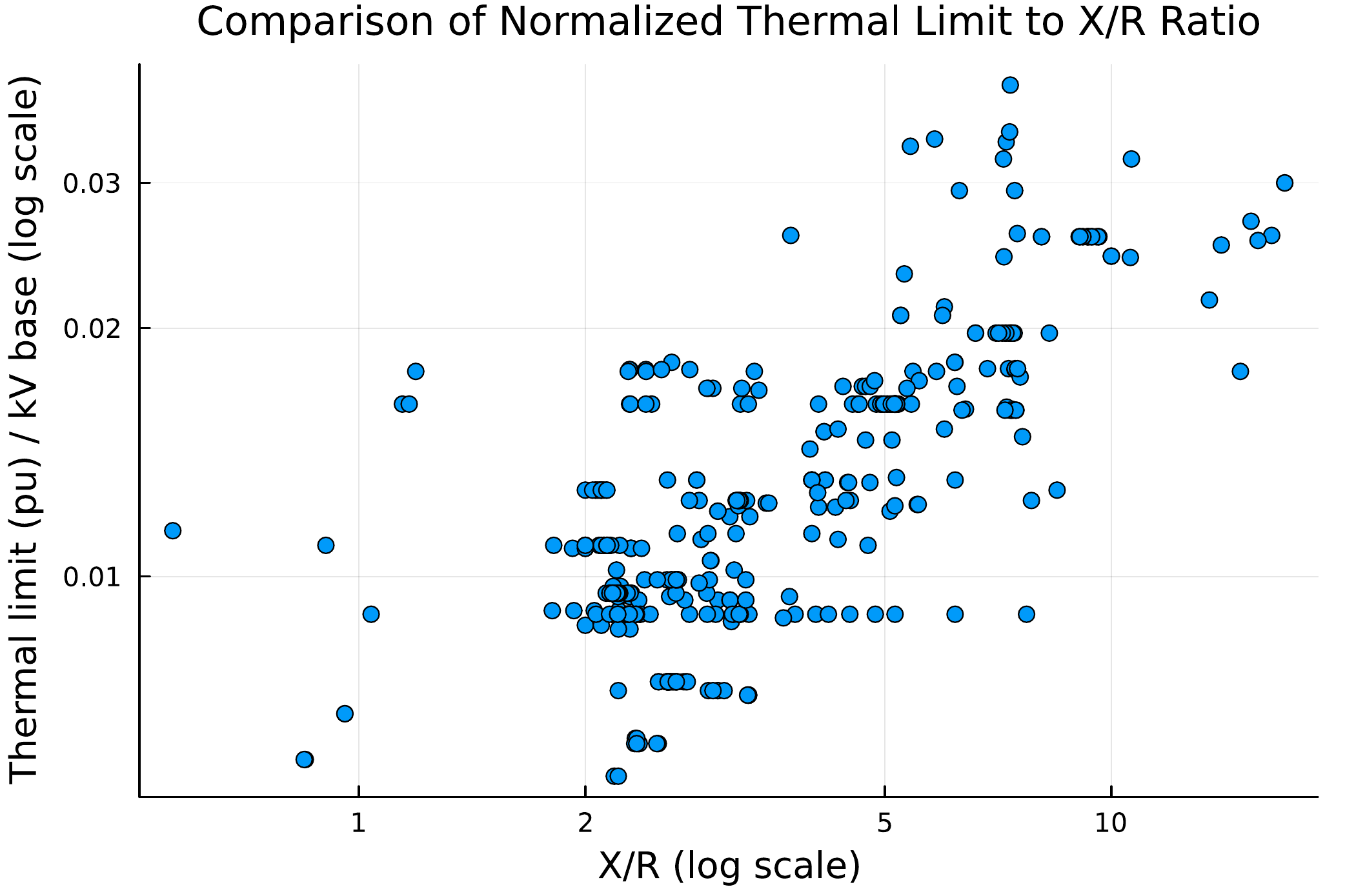}
  \caption{Polish network (left) and EIRGrid network (right) normalised thermal limits and X/R ratio correlation.}
  \label{fig:polish_irish_current_limit_xr_ratio}
\end{figure}

The statistical model proposed in \cite{Coffrin2014} derives thermal limits in power flow (MVA) rather than current flow (kA), which means that identical conductors can have different MVA limits based on their nominal voltage levels,  but converting back to kA they should have similar thermal limits. The model also utilizes the ratio of line reactance to resistance to de-correlate the line conductor type with line length. Figure~\ref{fig:polish_irish_current_limit_xr_ratio} indicates a clear correlation between $X/R$ ratio and the kA thermal limits.

\begin{figure}[H]
  \centering
  \includegraphics[width=0.48\columnwidth]{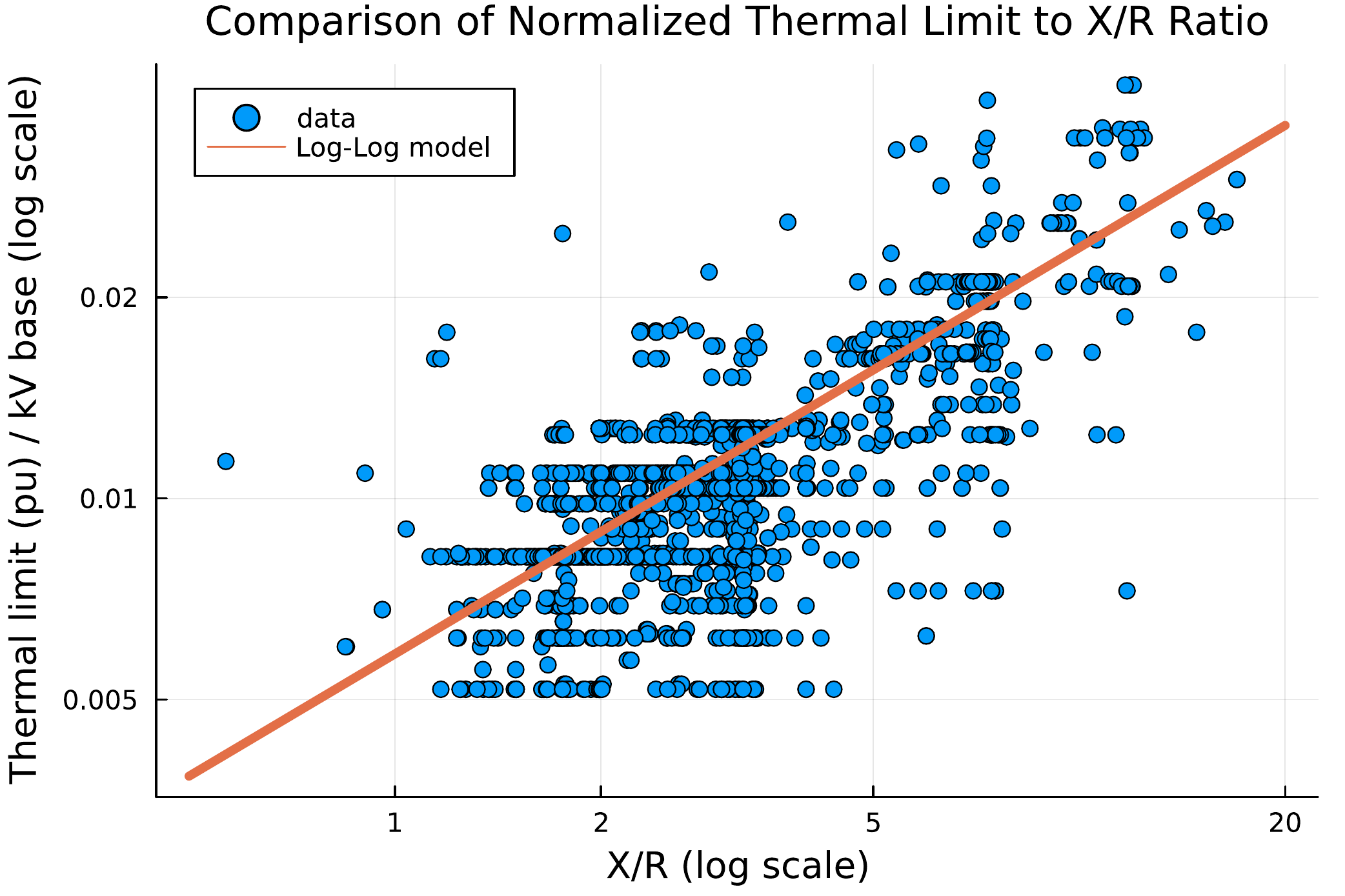}
  \caption{Linear log-log model for combined Polish and EIRGrid networks.}
  \label{fig:polish_irish_linear_loglog}
\end{figure}

We next combine the two datasets and fit the following linear regression model on the log-log scale
\begin{equation}
	\label{eq:linear_fit}
	y = e^a x^k
\end{equation}
for which we derive the following values for $a$ and $k$,
\begin{equation}
	a = -5.1407,  \quad 
	k = 0.6078.
\end{equation}
The regression fit model and its parameters are shown in Figure~\ref{fig:polish_irish_linear_loglog} which is a fairly crude but effective fit of the data for our purposes.
In \cite{Coffrin2014} the values for $a$ and $k$ are derived as $-5.0886$ and $0.4772$, respectively. Therefore the value for $a$ we derived is close to the one in \cite{Coffrin2014}, but the value of $k$ is slightly different, which could be due to having a more recent Irish dataset, inclusion of more points from the reference datasets and/or the filters applied to the data. 

The model \eqref{eq:linear_fit} is used to estimate per unit thermal limits $S^u$ for the lines that have all three parameters $X, R, v$ on both sides of the line as,
\begin{equation}
	\label{eq:statistical_model}
	S^u = v e^{-5.1407} \left( \frac{X}{R} \right) ^{0.6078}.
\end{equation}

\subsubsection{Upper Bound}
\label{sec:upper_bound}
The statistical model \eqref{eq:statistical_model} is only applicable when all three $X, R, v$ are available and both sides of the element have the same nominal voltage.  It is therefore not possible to apply this model to ideal lines with $R=0$ and transformers with different nominal voltage at either side. In these cases, the authors in \cite{Coffrin2014} propose calculating a reasonable theoretical upper bound for the thermal limit that is close to the theoretical maximum, but not to the extent that it deactivates the thermal limit constraints.  The upper bound accurately indicates the line's true throughput limitation and is calculated from \cite{Coffrin2014}
\begin{equation}
	(S^u_{ij})^2 = (v_i^u)^2 Y_{ij}^2 \big( (v_i^u)^2 + (v_j^u)^2 - 2 v_i^u v_j^u cos(\theta_{ij}^\Delta) \big)
\end{equation}
where the index $ij$ means the line/transformer between buses $i$ and $j$,  $v_i^u$ is the voltage magnitude upper bound at bus $i$, $Y_{ij}$ is the line/transformer's admittance, and $\theta_{ij}^\Delta$ is line/transformer's phase angle difference bound.  The thermal limit upper bound can be calculated given that reasonable voltage bounds are available throughout the network and with the assumption of $\theta_{ij}^\Delta = 15^{\circ}$.

\subsubsection{Validating the Statistical and Upper Bound Models}
\label{sec:validating_the_models}
\begin{figure}[H]
  \centering
  \includegraphics[width=0.48\columnwidth]{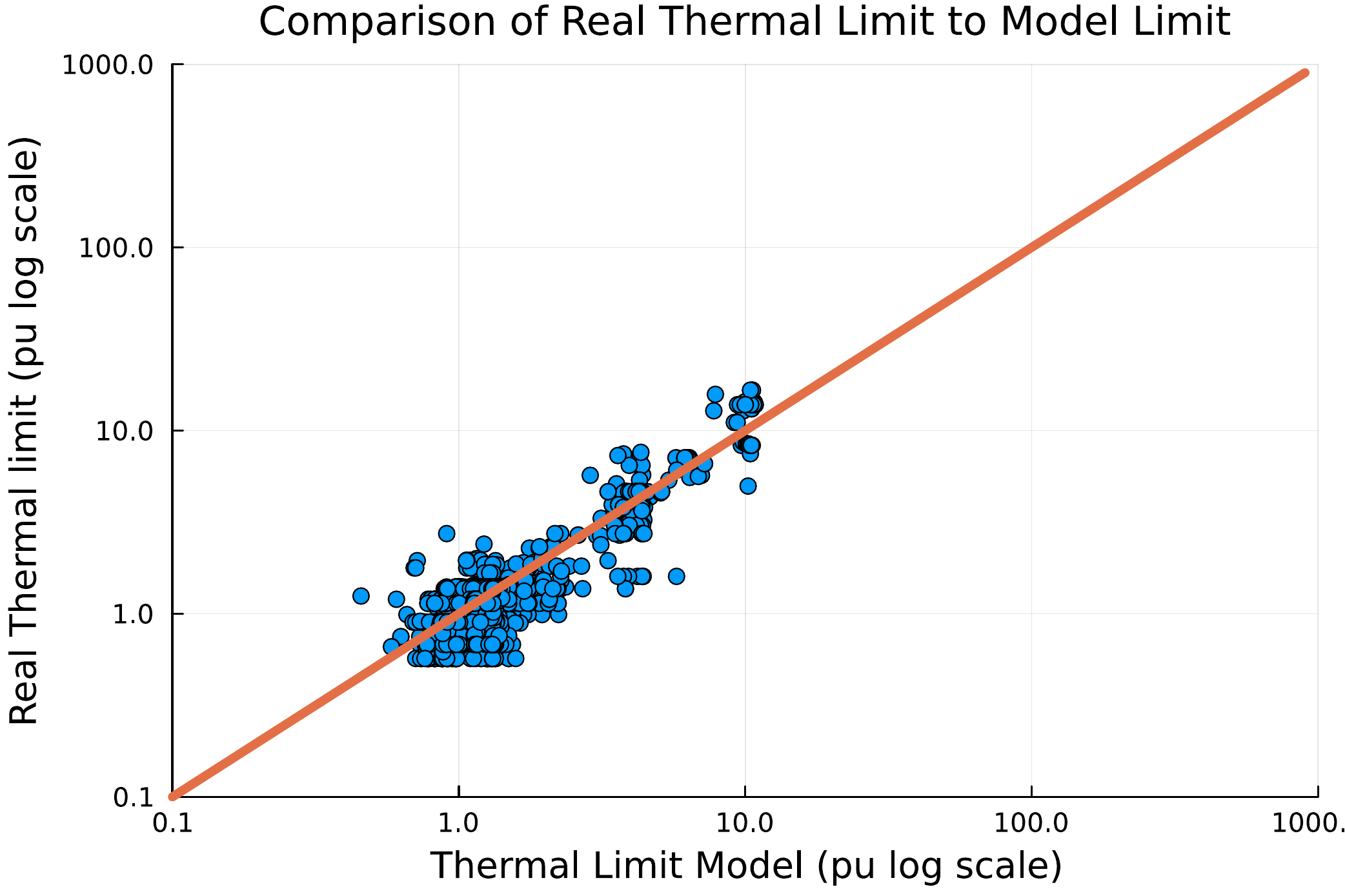}
    \includegraphics[width=0.48\columnwidth]{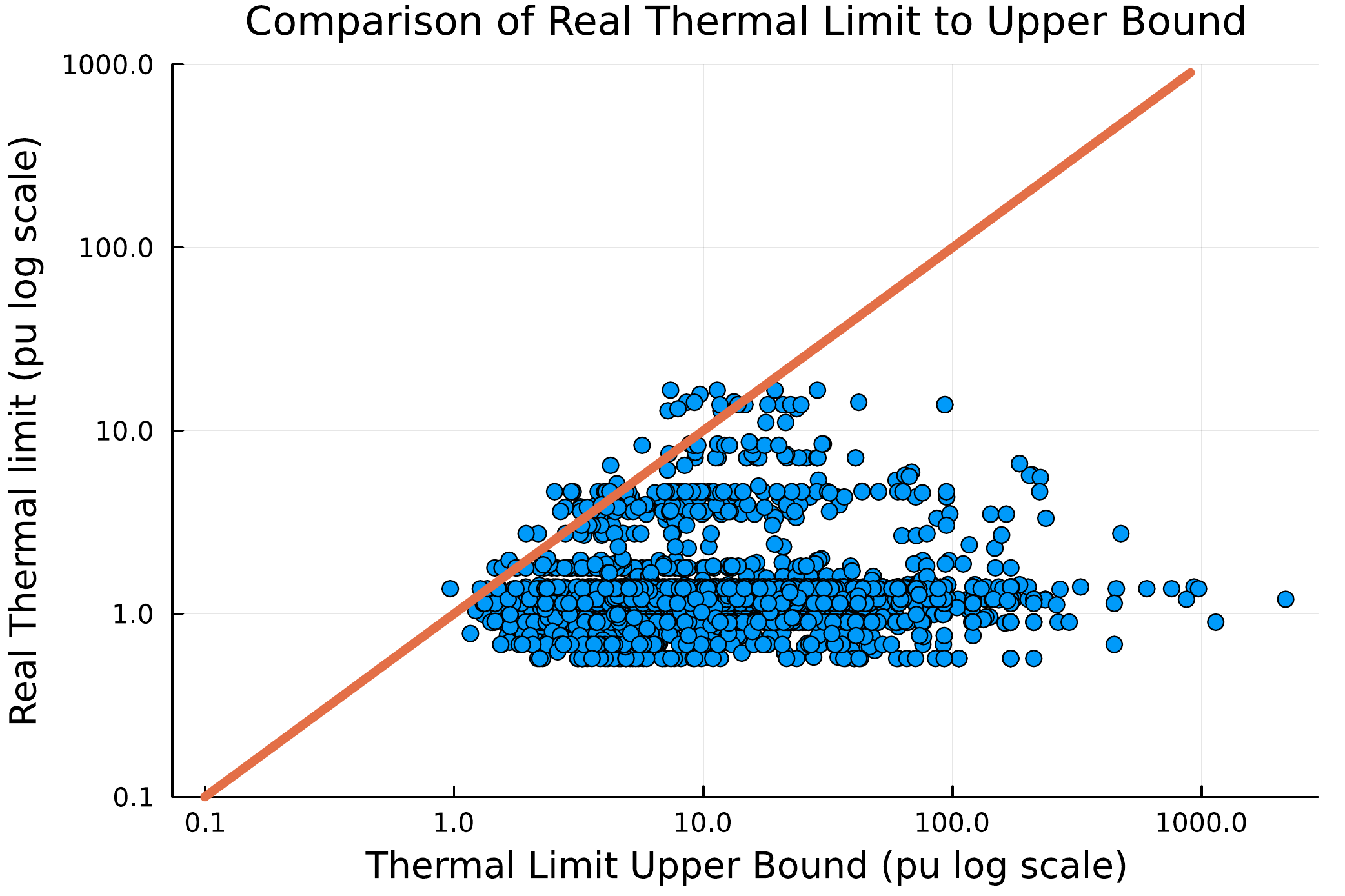}
  \caption{Comparison of the statistical model (left) and the thermal limit upper bound calculation (right) to the given MVA limits for Polish and EIRGrid networks.}
  \label{fig:polish_irish_validation}
\end{figure}

Figure~\ref{fig:polish_irish_validation} shows correlation of the proposed statistical and upper bound thermal limit values to the real thermal limits given in the reference datasets of Polish and EIRGrid transmission systems.  The figure on the statistical model shows a reasonable correlation of the model and the actual data as the data forms a stretch of points along the $x = y$ line which indicates perfect predictions.  The upper bound,  on the other hand, shows points below the $x = y$ line affirming the optimistic high bound on the thermal capacity by design,  while the large upper bound  limits indicate lines with small reactance values which represent ideal conductors.

\subsection{Reviewing the S-NEM2300 Data}
\label{sec:thermal_limits_NEM_data}
\begin{figure}[H]
  \centering
    \includegraphics[width=0.48\columnwidth]{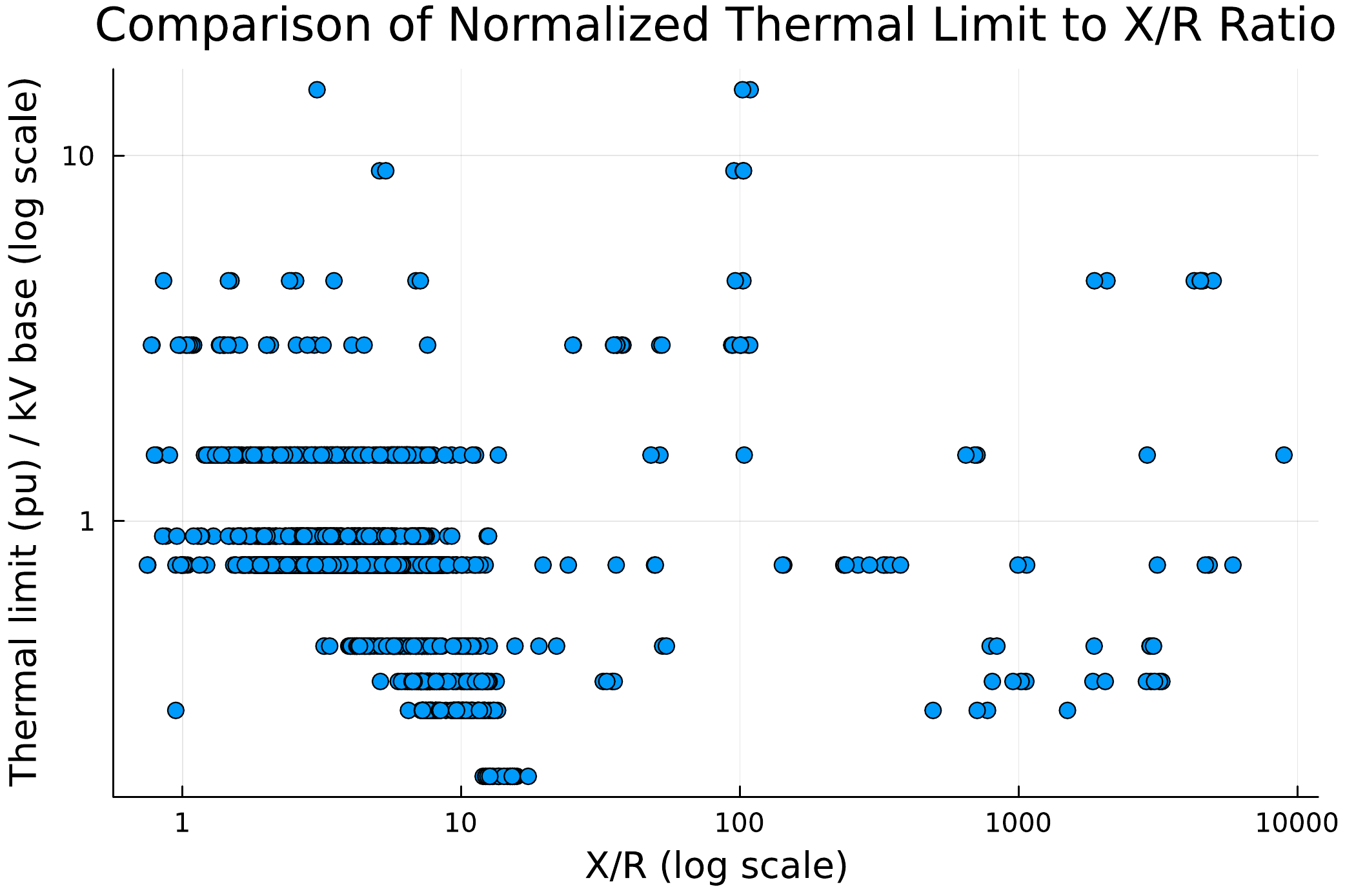}
  \caption{S-NEM2300 network normalised thermal limits and X/R ratio. Note that X/R ratio renders numerous infinity values corresponding to ideal lines which are not represented in this figure.}
  \label{fig:NEM_base_thermal_limit_kv_xr}
\end{figure}

Figure~\ref{fig:NEM_base_thermal_limit_kv_xr} shows the correlation between normalised thermal limits and the $X/R$ ratios in S-NEM2300 network. Note that all thermal limits are set to the nominal MVA of 100 MVA for all lines and transformers and only the $X/R$ values of non-ideal lines are shown.  From the figure it is evident the network data requires some modifications to show a correlation behaviour similar to reference datasets in Figure~\ref{fig:polish_irish_current_limit_xr_ratio}.

We perform two basic network reduction techniques to modify and prepare the network data for the thermal limit models application:
\begin{itemize}
	\item we remove the lines that represent ideal lines or circuit breakers, filtered as lines with
	\begin{itemize}
		\item small impedance ($<0.01$ p.u.) and small shunt admittance ($<0.01$ p.u.),  since in the synthetic data there exist constant parameter (CP) and decoupling lines with negligible impedance but non-negligible shunt admittance values, or,
		\item high $X/R$ ratio ($>100$) and small shunt admittance ($<0.01$ p.u.);
	\end{itemize}
	\item we join the cascading lines connected through degree 2 buses to ensure a uniform thermal limit over linear paths in the network.
\end{itemize}

\begin{figure}[H]
  \centering
    \includegraphics[width=0.48\columnwidth]{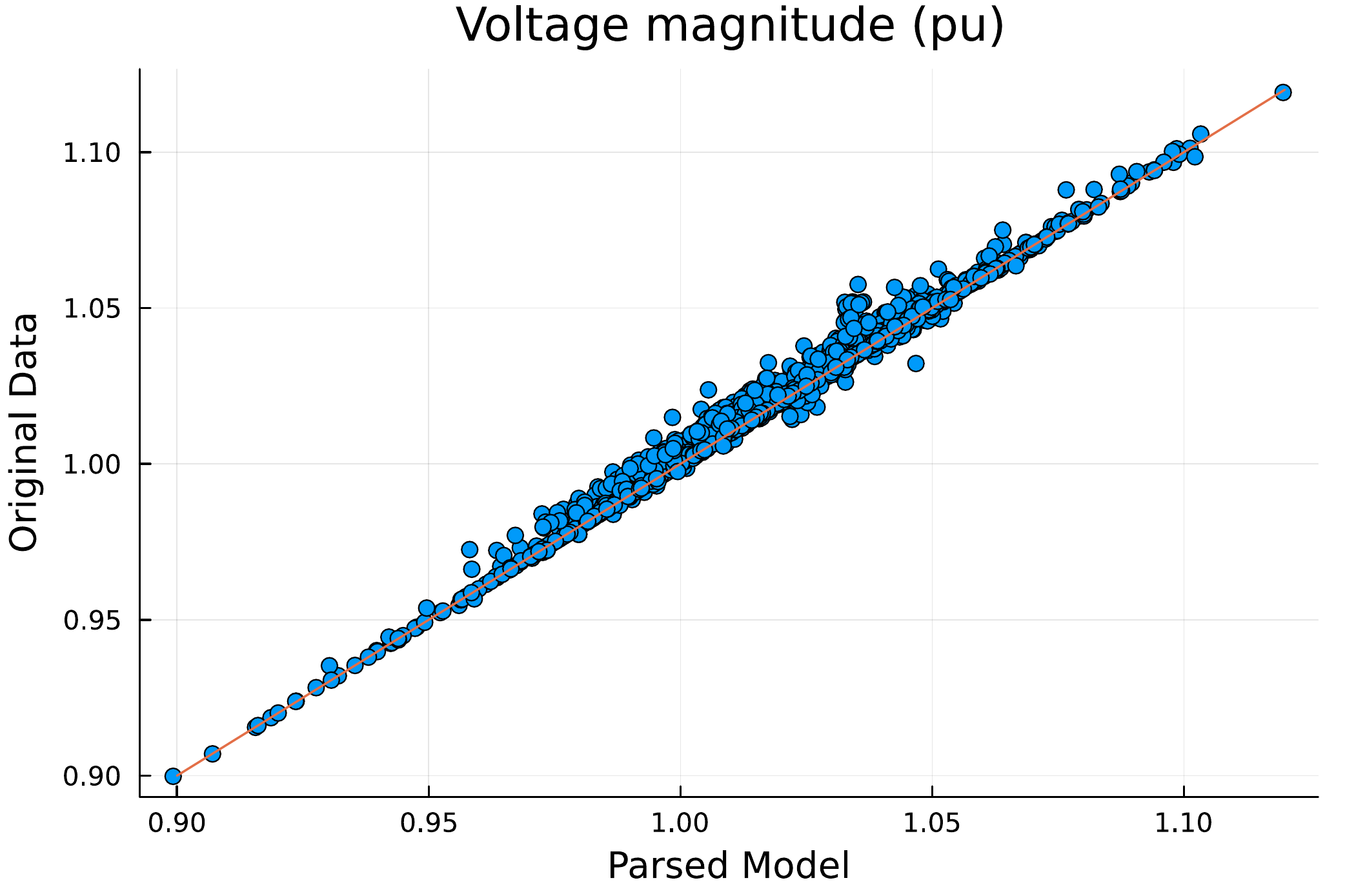}
    \includegraphics[width=0.48\columnwidth]{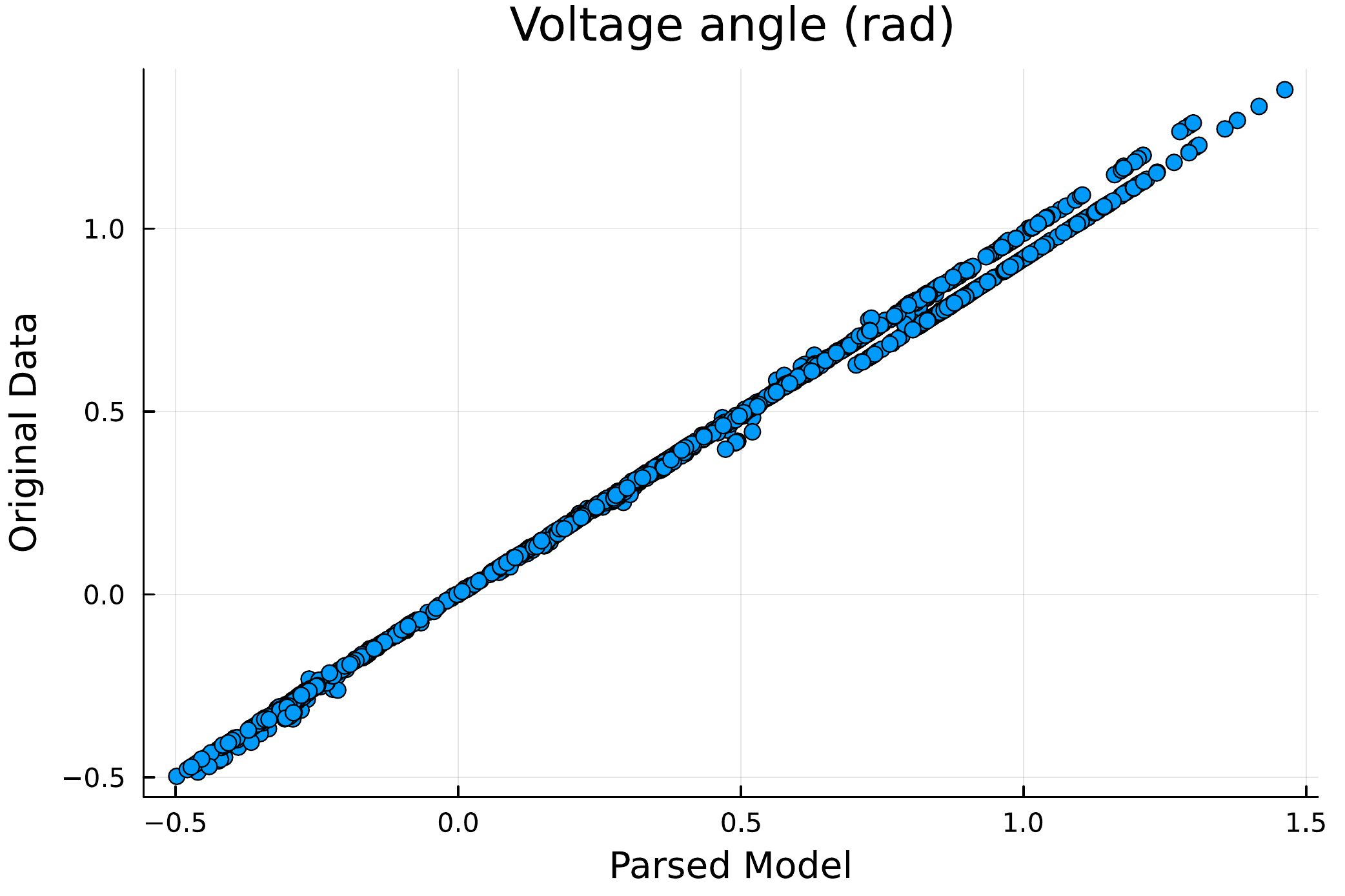}
  \caption{Validation of S-NEM2000 after network reduction.}
  \label{fig:validation_reduced_vm_va}
\end{figure}

Figure~\ref{fig:validation_reduced_vm_va} shows how the reduced network bus voltage magnitudes and phase angles results still closely match the original network results. In fact,  we used this power flow validation to specify and visualise a small variation of voltage magnitude and angle w.r.t the S-NEM2300 results in order to determine the empirical thresholds for impedance,  shunt admittance, and $X/R$ ratio values of the lines to be removed.

\subsection{Adding Thermal Limits and Validations}
\label{sec:adding_thermal_limits}
\begin{figure}[H]
  \centering
  \includegraphics[width=0.48\columnwidth]{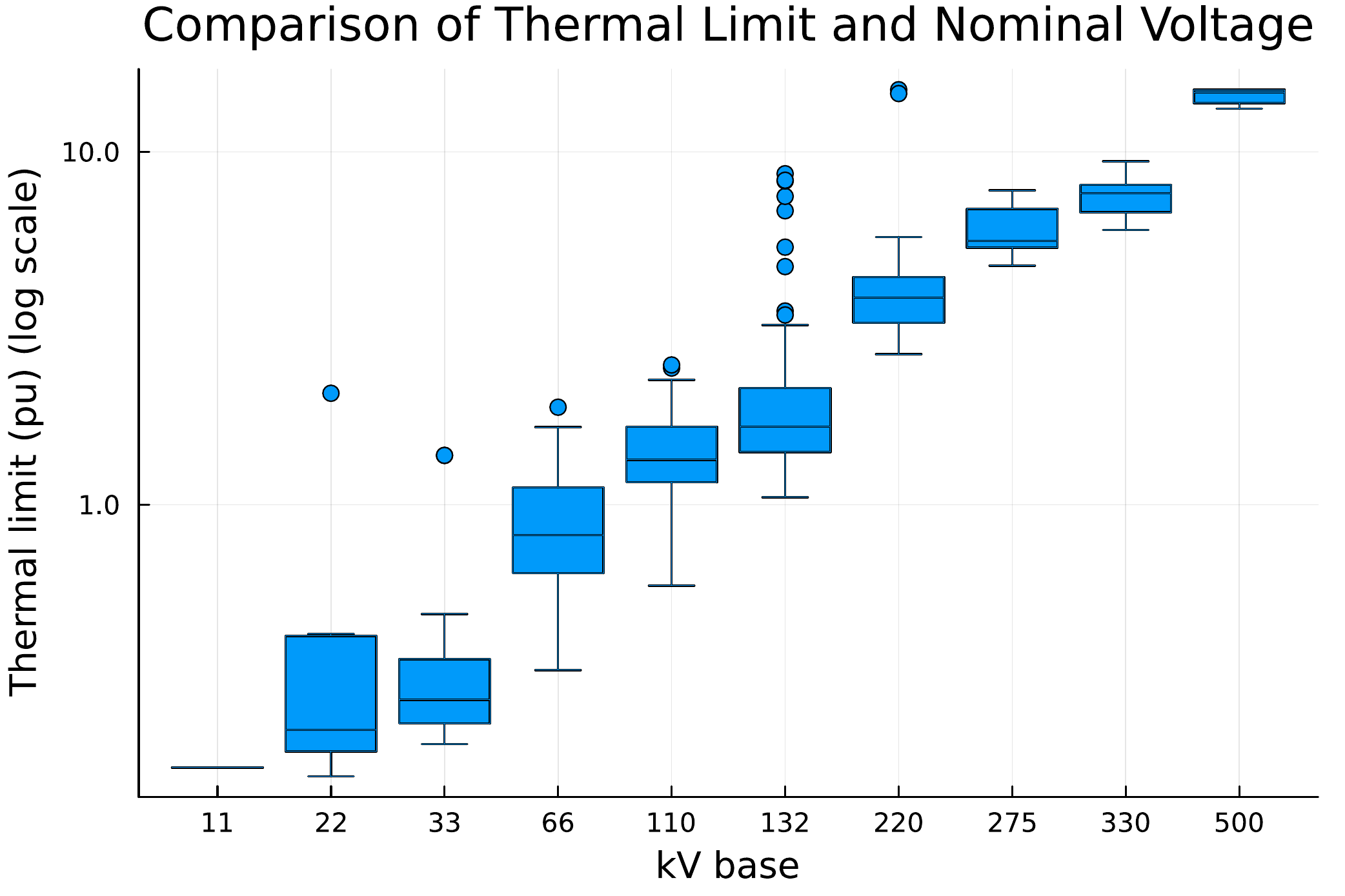}
    \includegraphics[width=0.48\columnwidth]{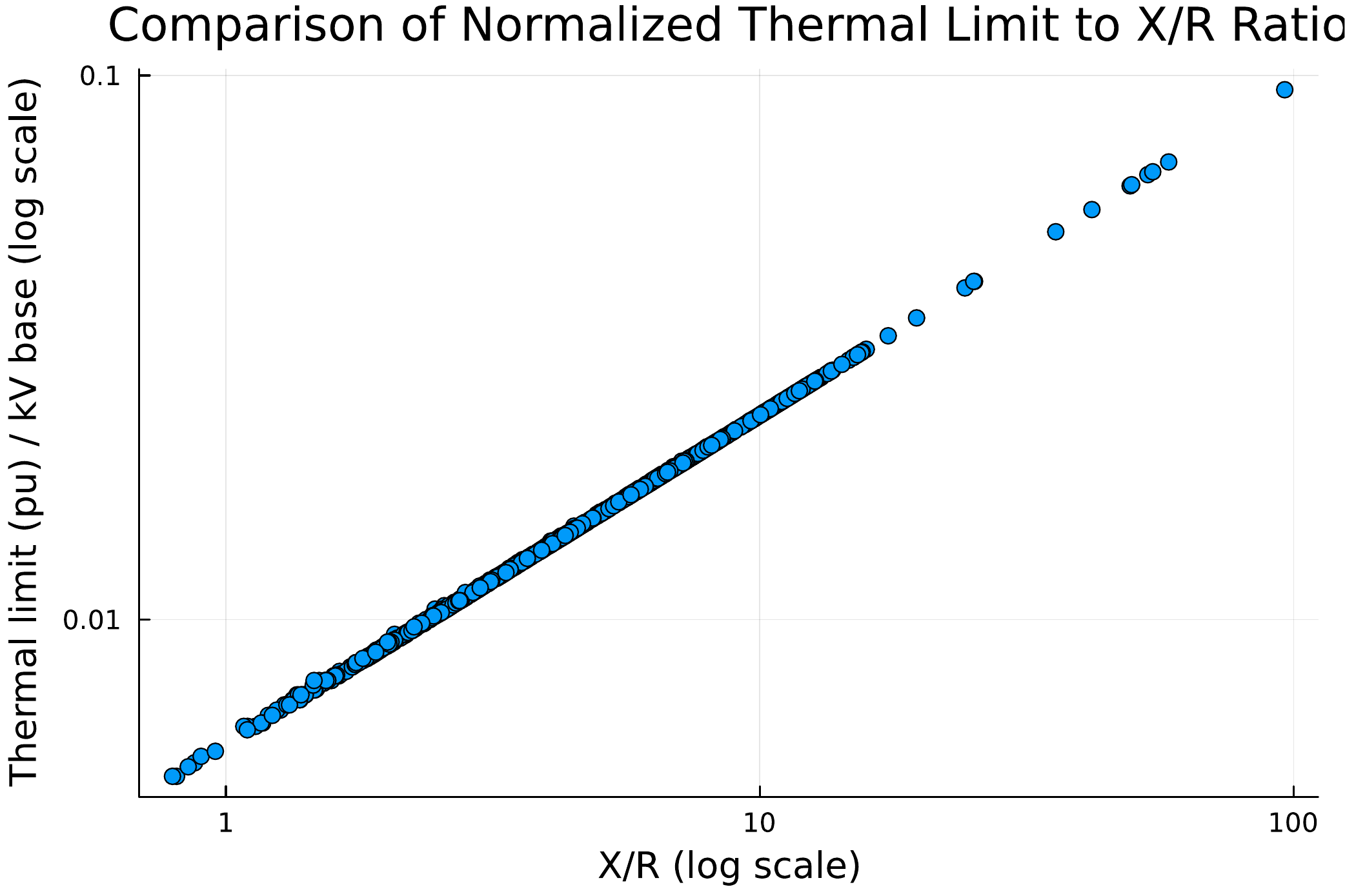}
  \caption{S-NEM2000 thermal limit vs voltage levels (left) and normalised thermal limits and X/R ratio (right). }
  \label{fig:NEM_not_fixed_thermal_limit_kv_xr}
\end{figure}

In the next step we add line and transformer thermal limit using the statistical and upper bound thermal limit models where applicable. 
We also set lower and upper bounds of 30 and 1500 MVA for the thermal limits on transformers according to NEM equipment ratings \cite{AEMO_ratings}.
Figure~\ref{fig:NEM_not_fixed_thermal_limit_kv_xr} shows line thermal limits across different voltage levels and correlation of normalised thermal limits and $X/R$ ratio.  The left figure shows that the lines with higher nominal voltage have higher thermal limits which is a true behaviour of real systems. The perfect prediction of line thermal limits in the right figure indicates that no ideal lines exists in the data any more.  Note that the transformer thermal limits are derived from the upper bound model and are not represented in Figure~\ref{fig:NEM_not_fixed_thermal_limit_kv_xr}.

\begin{figure}[htbp]
  \centering
    \includegraphics[width=0.43\columnwidth]{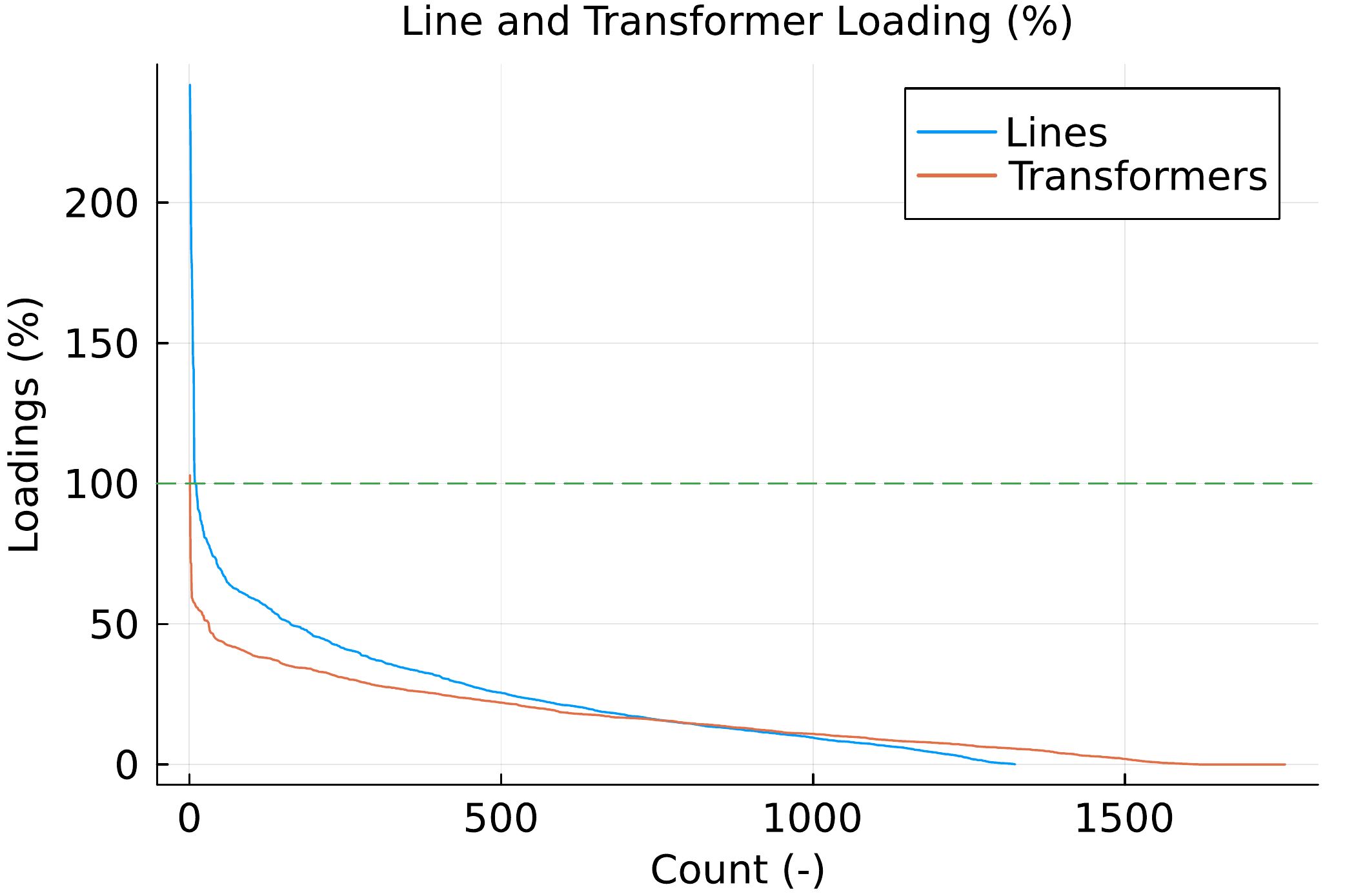}
    \includegraphics[width=0.43\columnwidth]{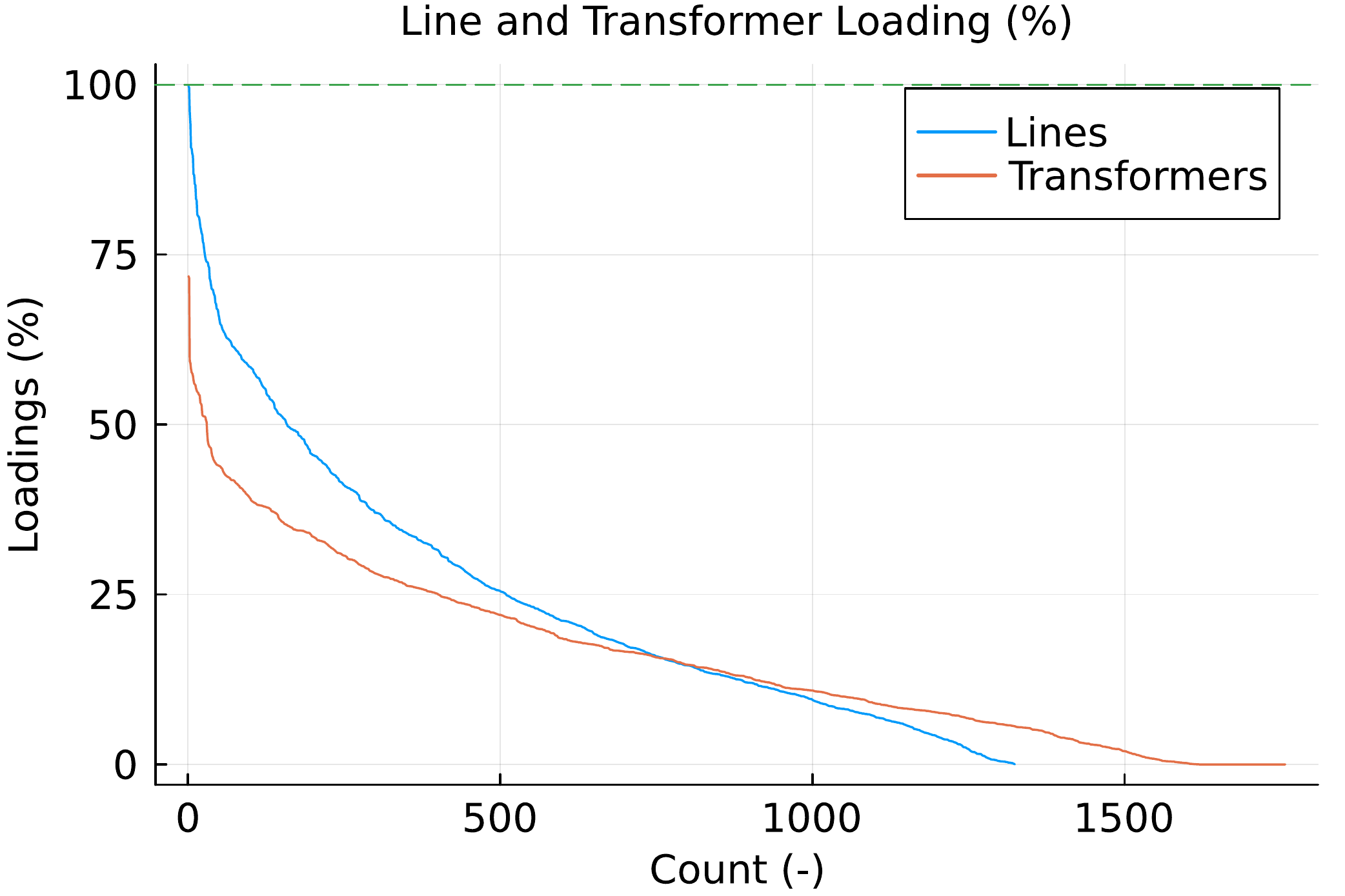}
  \caption{S-NEM2000 line and transformer loading: before modification (left), after modification(right). }
  \label{fig:NEM_loading}
\end{figure}

\begin{figure}[htbp]
  \centering
  \includegraphics[width=0.43\columnwidth]{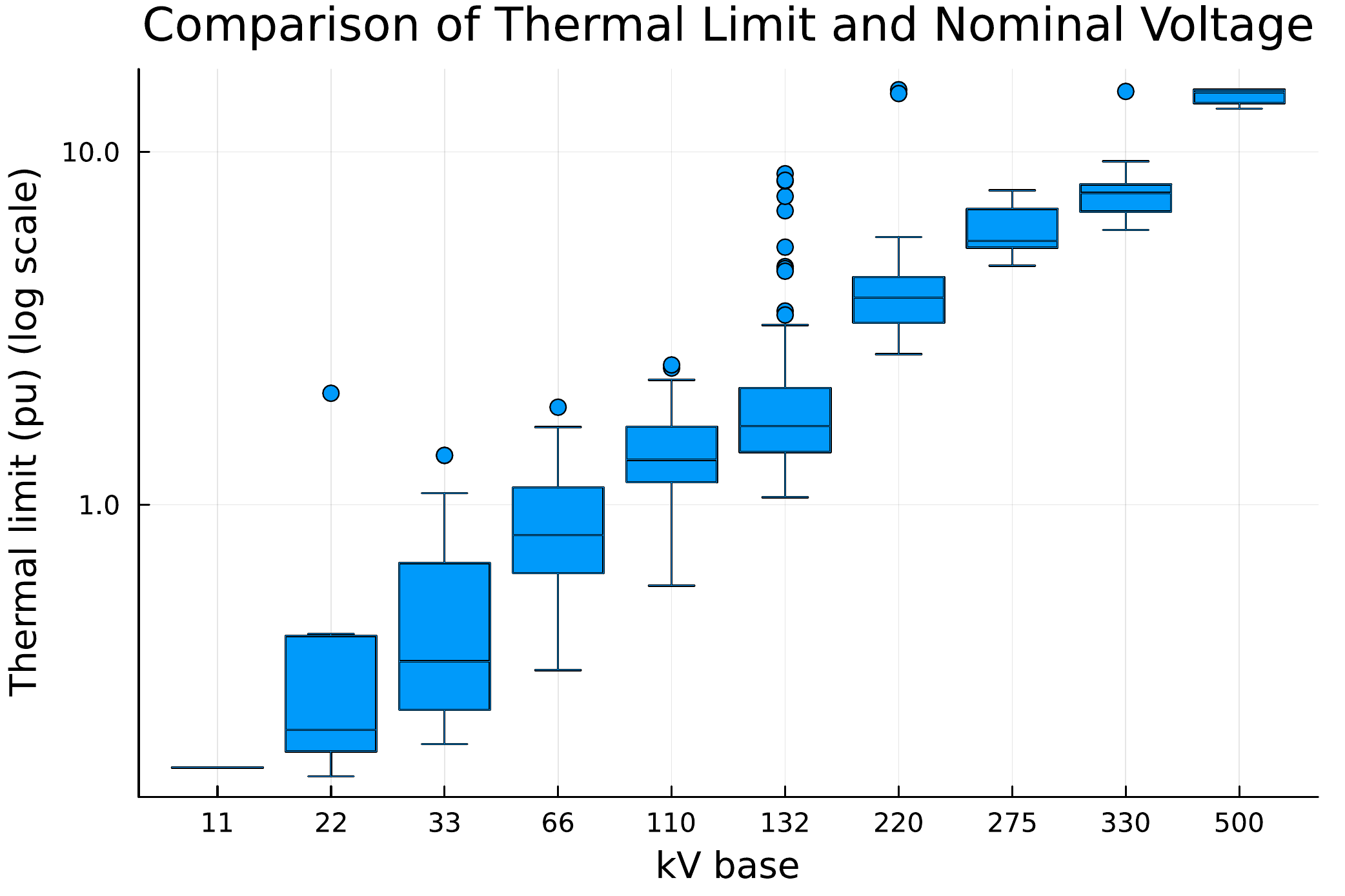}
    \includegraphics[width=0.43\columnwidth]{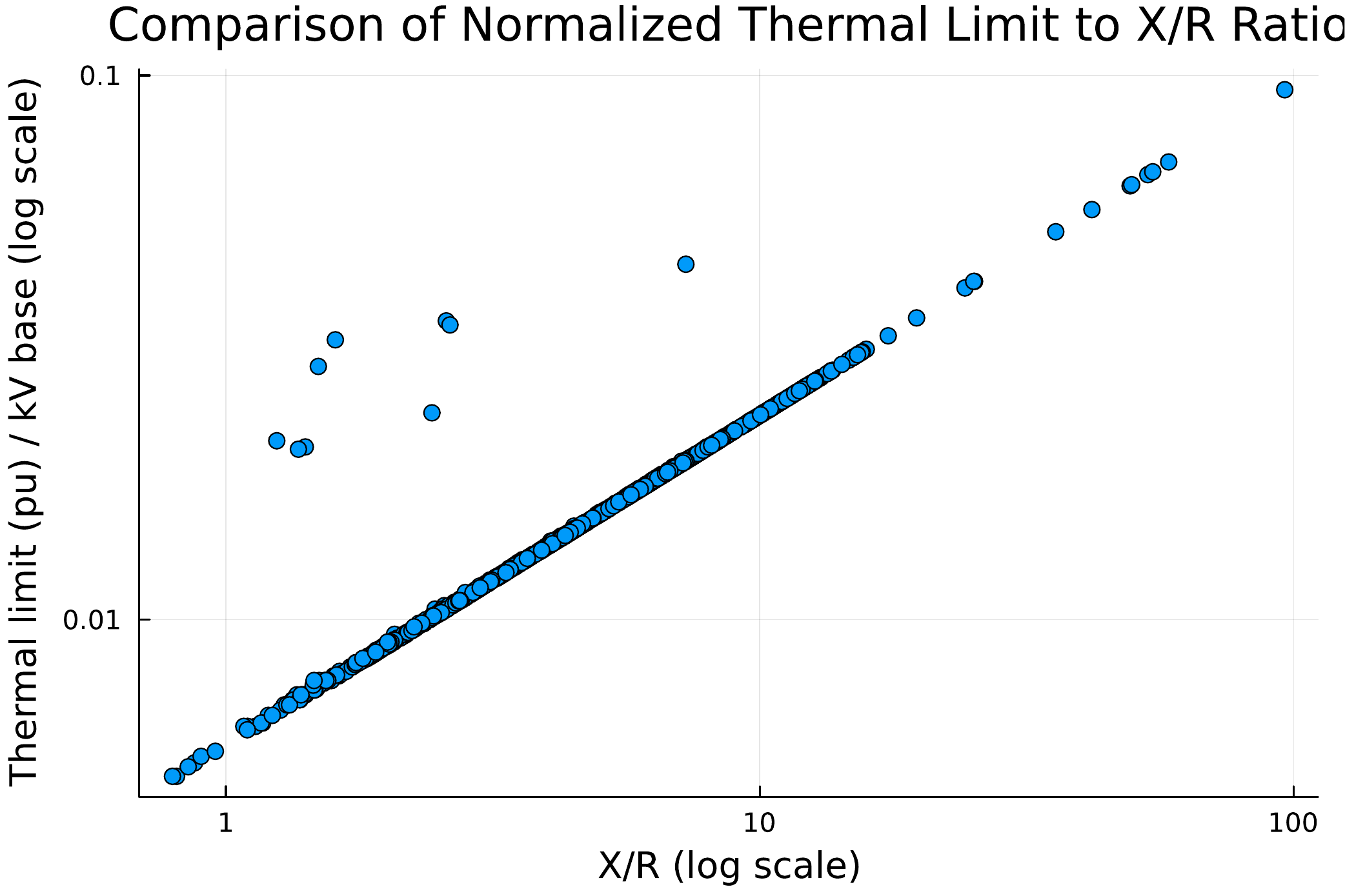}
  \caption{S-NEM2000 modified thermal limit vs voltage levels (left) and normalised thermal limits and X/R ratio (right).}
  \label{fig:NEM_fixed_thermal_limit_kv_xr}
\end{figure}

The estimated thermal limits, though seem to have been predicted reasonably well as seen in Figure~\ref{fig:NEM_not_fixed_thermal_limit_kv_xr}, should be tested against power flow results of the modified network data. Figure~\ref{fig:NEM_loading} on the left shows line and transformer loading w.r.t their thermal limits. We note that with these thermal limits a total of 9 line flow violations occur in steady state conditions.
Since the number and magnitude of violations are manageable, in the next step we increase the thermal limits of the violating lines appropriately to ensure that the network is secure under normal operation. Figure~\ref{fig:NEM_loading} on the right shows line and transformer loading w.r.t their thermal ratings while Figure~\ref{fig:NEM_fixed_thermal_limit_kv_xr} shows how the thermal limits have changed with respect to Figure~\ref{fig:NEM_not_fixed_thermal_limit_kv_xr}.


\subsection{S-NEM2000 Network Model Components}
\label{sec:derived_syhtnetic_nem_data}
After cleaning the S-NEM2300 network data, reducing the lines, and fixing some minor issues with generators, the network data components count is listed in Table~\ref{table_derived_nem_elements}, which compared to the ones reported in Table~\ref{table_original_nem_elements}, the main differences are the number of buses and lines due to removing degree-2 buses and joining the connecting lines.  Note that the three-winding transformers are modelled as three two-winding transformers in the released version of the S-NEM2000 network data.

\begin{table}[tbh]
  \centering 
   \caption{S-NEM2000 network number of components.} 
   \label{table_derived_nem_elements}
    \begin{tabular}{l  r c }
    Component & Count \\
    \hline
    Buses & 2000 \\
    Lines & 1324 \\
    Two-winding transformers & 1418 \\
    Three-winding transformers & 113 \\
    Generators & 265\\
    Shunts & 301 \\
    Loads * & 1702 \\
\hline
\end{tabular}
\\
* Some generators have negative injection. \\
They are divided into a generator \\
with zero injection, and a load.
\end{table}


\section{Generation Models}
\label{sec:generation_models}
The canonical objective in OPF is generation cost minimization. Such an application, and related ones such as market clearing and unit commitment, depend on detailed generator specifications such as fuel type, active and reactive power capability curves, ramp up/down rates,  minimum on/off times, startup/shutdown costs, and fuel efficiency. To assign the appropriate properties to generators in the synthetic model, we first notice that most generator properties are directly linked to its mechanical design, which in turn is determined by the fuel type. We therefore develop data driven models to assign generator fuel categories, and then, using publicly available data sets we derive generator properties.

To assign generator fuel categories and properties we use publicly available generation datasets provided by AEMO \cite{AEMO2017, AEMO2014} and GeoScience Australia (GA) \cite{GA2021}. We categorize the generators fuel type so that the synthetic network model is a reasonable representative of the actual NEM generation fuel mix.

\subsection{Reviewing the Reference Data}
\label{sec:gen_reference_data}
The reference data we use to assign synthetic generator types is GA major power stations list \cite{GA2021} which contains generators' fuel type, registered capacity, voltage level, status, and more.  
%
We only consider operational generators and the ones with a single fuel source. The reference dataset \cite{GA2021} classifies generation fuel categories into black/brown coal,  gas, distillate, biomass, biogas, water, solar, wind,  battery, fuel oil and coal seam methane.  

\begin{figure}[H]
  \centering
  \includegraphics[width=0.7\columnwidth]{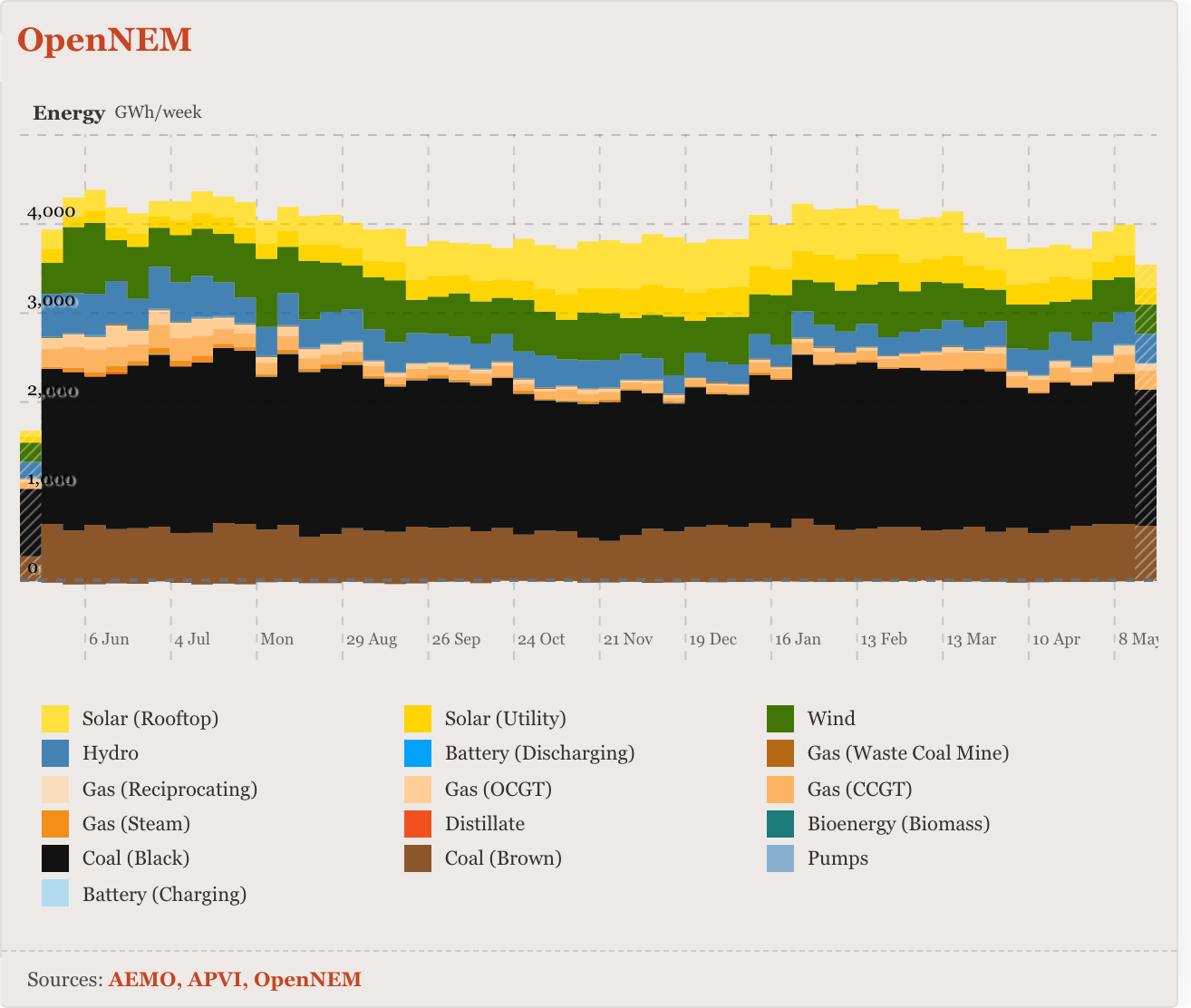}
  \caption{NEM generation fuel mix for the period of May 2022 to May 2023 \cite{OpenNem_license}.  }
  \label{fig:OpenNEM_fuel_mix}
\end{figure}

\begin{figure}[H]
  \centering
  \includegraphics[width=0.9\columnwidth]{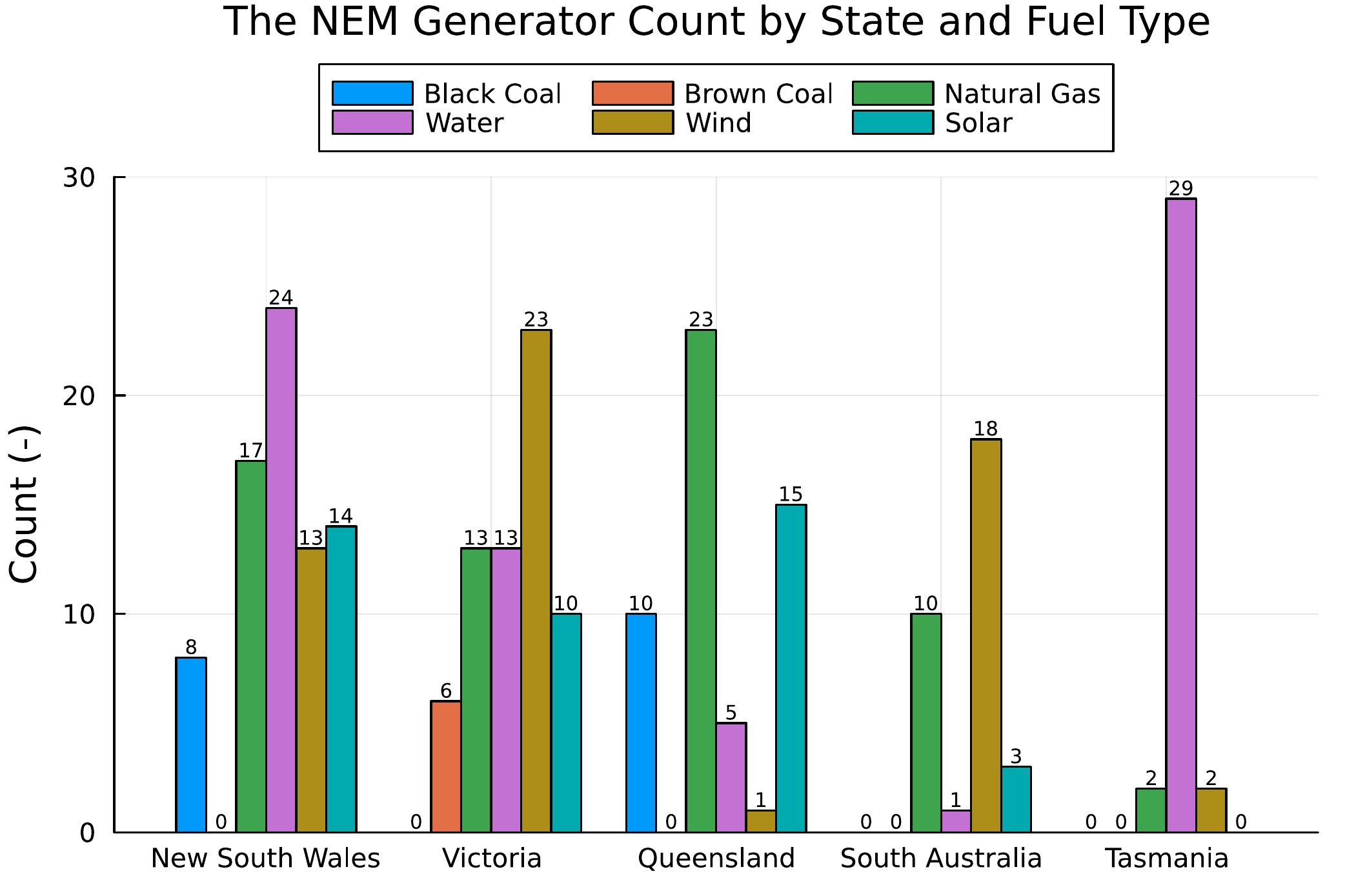}
  \caption{Generators count categorized by state and fuel type in the NEM.}
  \label{fig:gen_ref_count_by_fueltype_by_state}
\end{figure}

\begin{figure}[H]
  \centering
    \includegraphics[width=0.9\columnwidth]{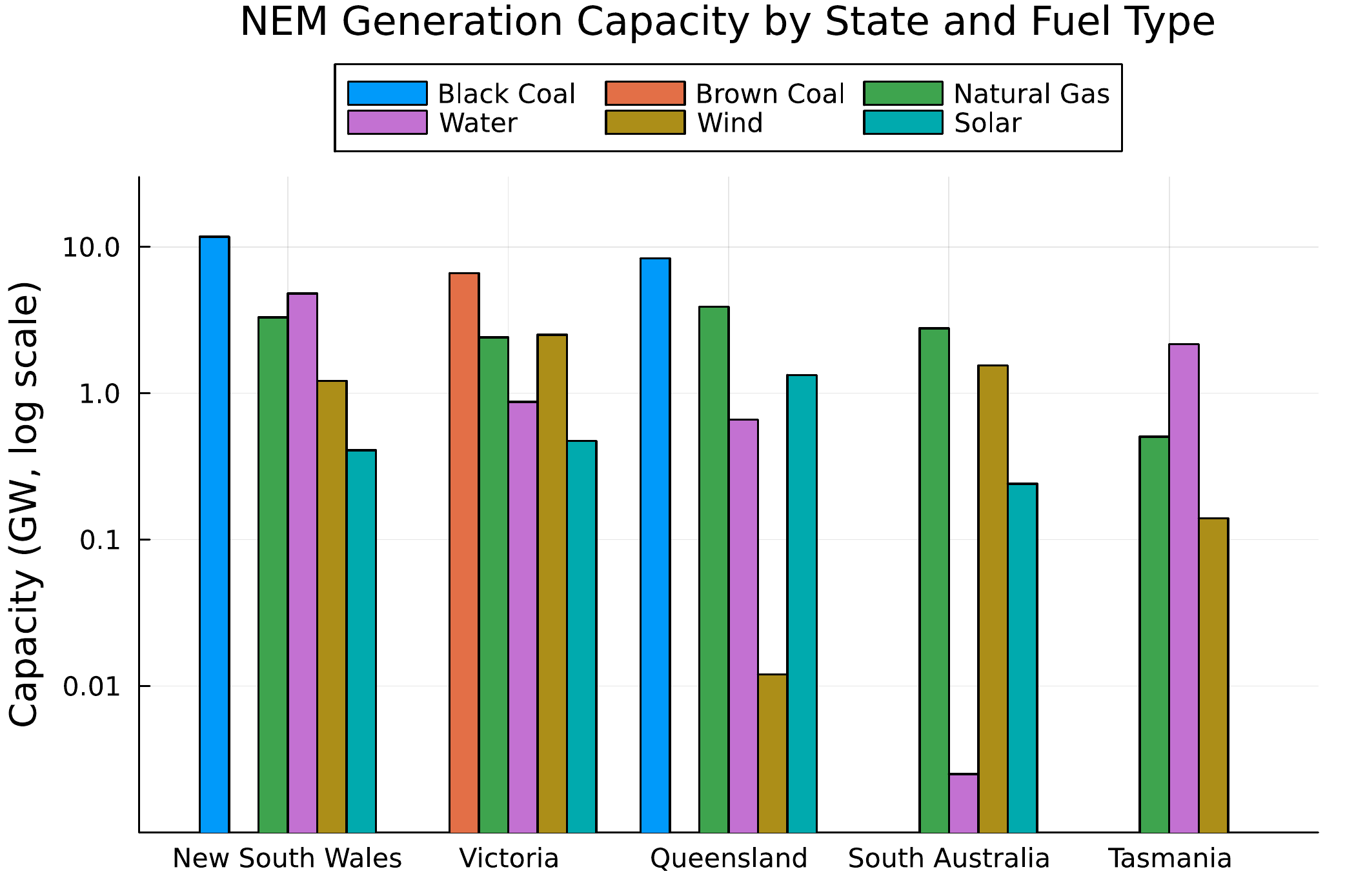}
  \caption{Generation capacity categorized by state and fuel type in the NEM.}
  \label{fig:gen_ref_capacity_by_fueltype_by_state}
\end{figure}

We first ignore generators using coal seam methane and fuel oil as there are only a few of them with very limited capacity.  Next,  inspired by Figure~\ref{fig:OpenNEM_fuel_mix} that shows generation fuel mix of the NEM for the past year \cite{OpenNem_license}, we ignore distillate, biomass, biogas and battery power stations as well.
Figures~\ref{fig:gen_ref_count_by_fueltype_by_state} and \ref{fig:gen_ref_capacity_by_fueltype_by_state} show the count and capacity of generators throughout the NEM categorized by fuel type and state.  For the S-NEM2000 model we target generation fuel type allocation to represent the generation fuel mix as shown in these two figures.

Upon determining generation fuel categories,  we assign generator properties and fuel cost data from Fuel and Technology Cost Review \cite{AEMO2014} which contains parameters such as start up and shut down costs, no-load, fixed and variable operation costs, minimum up/down times, ramp up/down rates, and more operational properties. 

%
%
%

\subsection{Reviewing the S-NEM2000 Data}
\label{sec:gen_nem_data}
The S-NEM2300 model classifies generators into synchronous generators and network sources, and does not take into account any fuel category and cost data. However, the authors in \cite{ArranoVargas2022} state that synchronous generations include fossil-fuel generators as well as hyrdo and wind powered generators whereas the network source category consists of `other' types of generators. We base our classification on these assumptions and choose the network sources to only include solar farms.  It is also worth noting that the reactive power compensators such as synchronous var generators and synchronous condensers can be part of the network sources as well, particularly by noting that some of the networks sources only provide reactive power and no active power to the network.

\paragraph{Remark:}
The S-NEM2300 network data models the Basslink HVDC interconnector between the mainland NEM and Tasmania as two network sources connected to each sub-network. However, since both network sources are indexed and placed as Tasmania elements, we assign the generator (and the generator bus), representing Tasmania in the mainland, to the Victoria area.

\begin{figure}[H]
  \centering
  \includegraphics[width=0.47\columnwidth]{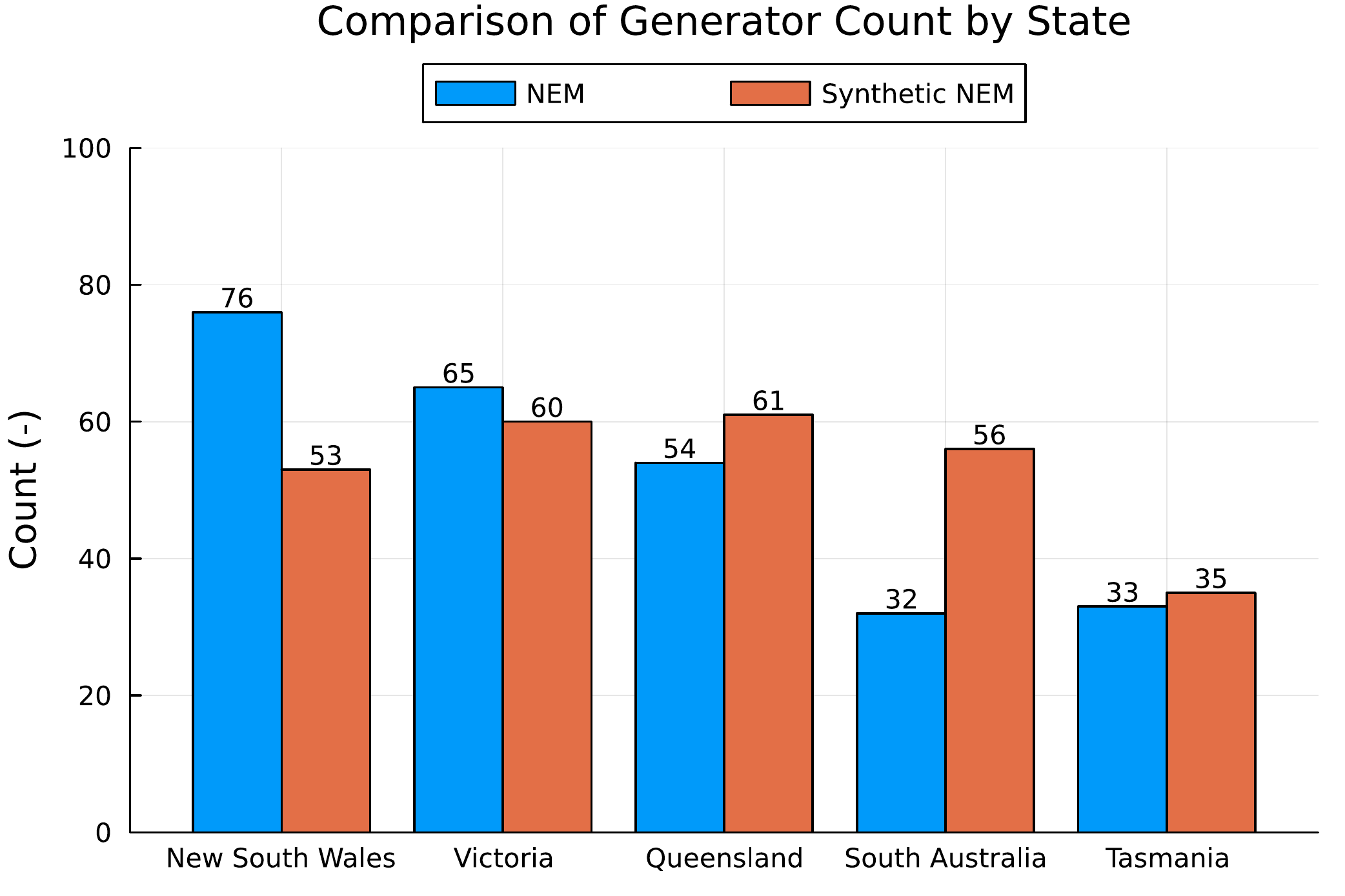}
  \includegraphics[width=0.47\columnwidth]{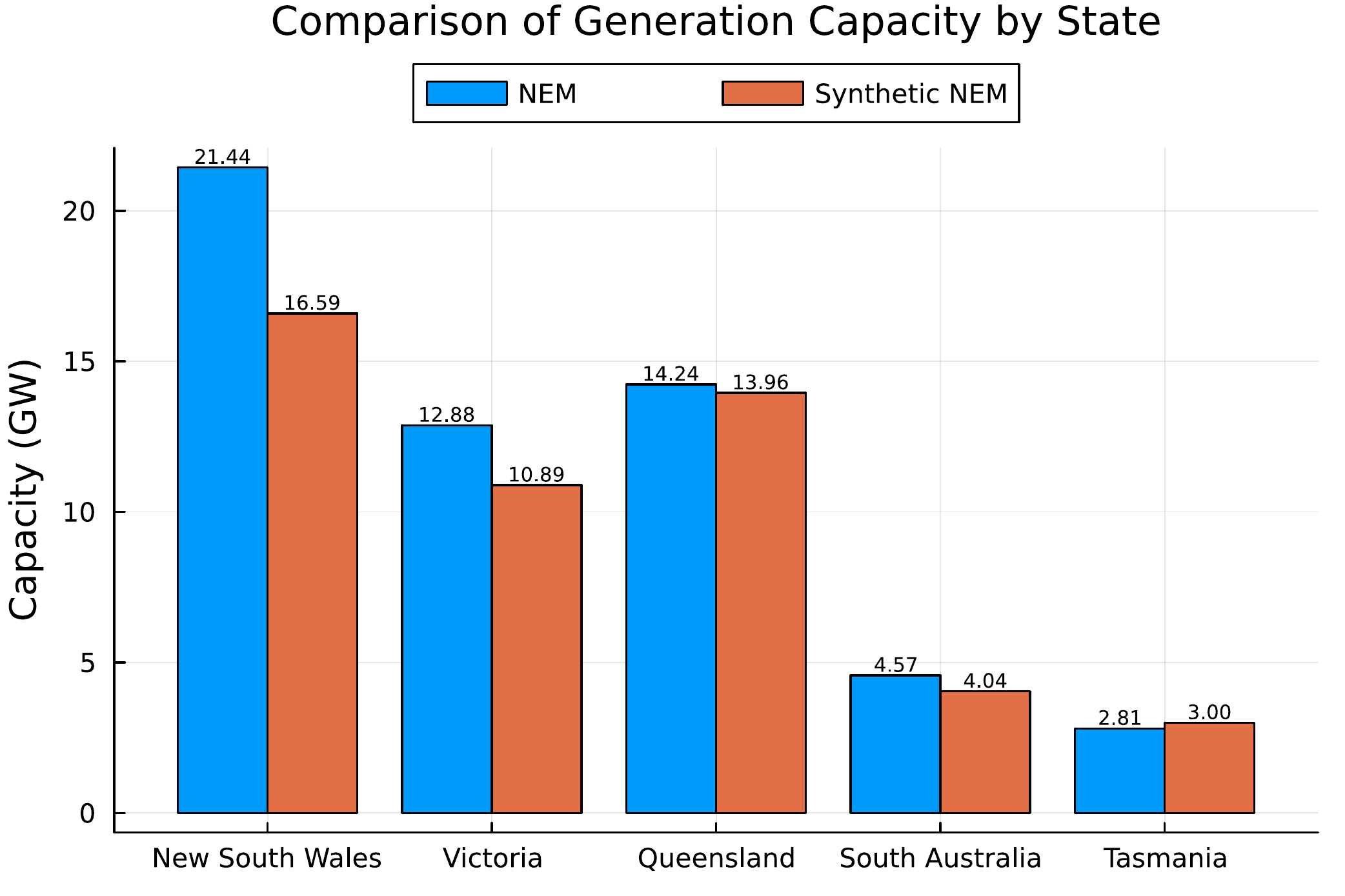}
  \caption{Comparison of generators in the NEM and Synthetic NEM per state with respect to: total generators count (left) and total generation capacity (right).}
  \label{fig:gen_compare_capacity_count_by_fueltype_by_state}
\end{figure}

Figure~\ref{fig:gen_compare_capacity_count_by_fueltype_by_state} shows how the total count and capacity of generators across different states are compared between the NEM and the S-NEM2300.

\subsection{Generation Fuel Category Classification}
In effort to classify synthetic generators by fuel category, a data-driven classification approach was taken. Within this approach the registered capacity of generators from GA major power stations list \cite{GA2021} were used to train an Sklearn Random Forest classification model \cite{sklearn_api, breiman2001random} to predict the fuel type of a given generation.  Given that in the reference dataset some fuel types are overrepresented compared to others,  we used oversampling of the under-represented classes by generating samples based on nearest neighbours using Synthetic Minority Over-sampling Technique (SMOTE) \cite{SMOTE}.  Using SMOTE to balance the classes before classifying, the trained classification model was then used with the maximum generation capacity of synthetic generators within the S-NEM2000 dataset to predict the fuel type of the generators in each state separately.
%

To stay true to the S-NEM2300 model, the classifier uses the synthetic synchronous generators list to identify generators with black/brown coal, natural gas, hydro and wind fuel categories and uses the network source list to assign solar generation units\footnote{
An important point that is missing from many datasets is the identification of synchronous condensers and other var compensator units. The S-NEM2300 model does not mention var compensators as part of the generator categories, but close observation of the data reveals that the network sources contains several generation units with zero or very small active and substantial reactive power contribution to the network.  One could assume these units as var compensator units and assign the reminder as solar generation units. 
}.
Nonetheless, we pre-allocate the the network source representing Basslink in NEM mainland as a hydro generator, due to hydro being the main fuel category in Tasmania, and similarly the network source representing Basskink in Tasmania as a natural gas fueled generator, due to natural gas being the only fossil-fuel present in Tasmania fuel categories\footnote{These two network sources can be removed in future work when HVDC lines are added to the network.}.

\begin{figure}[H]
  \centering
  \includegraphics[width=0.9\columnwidth]{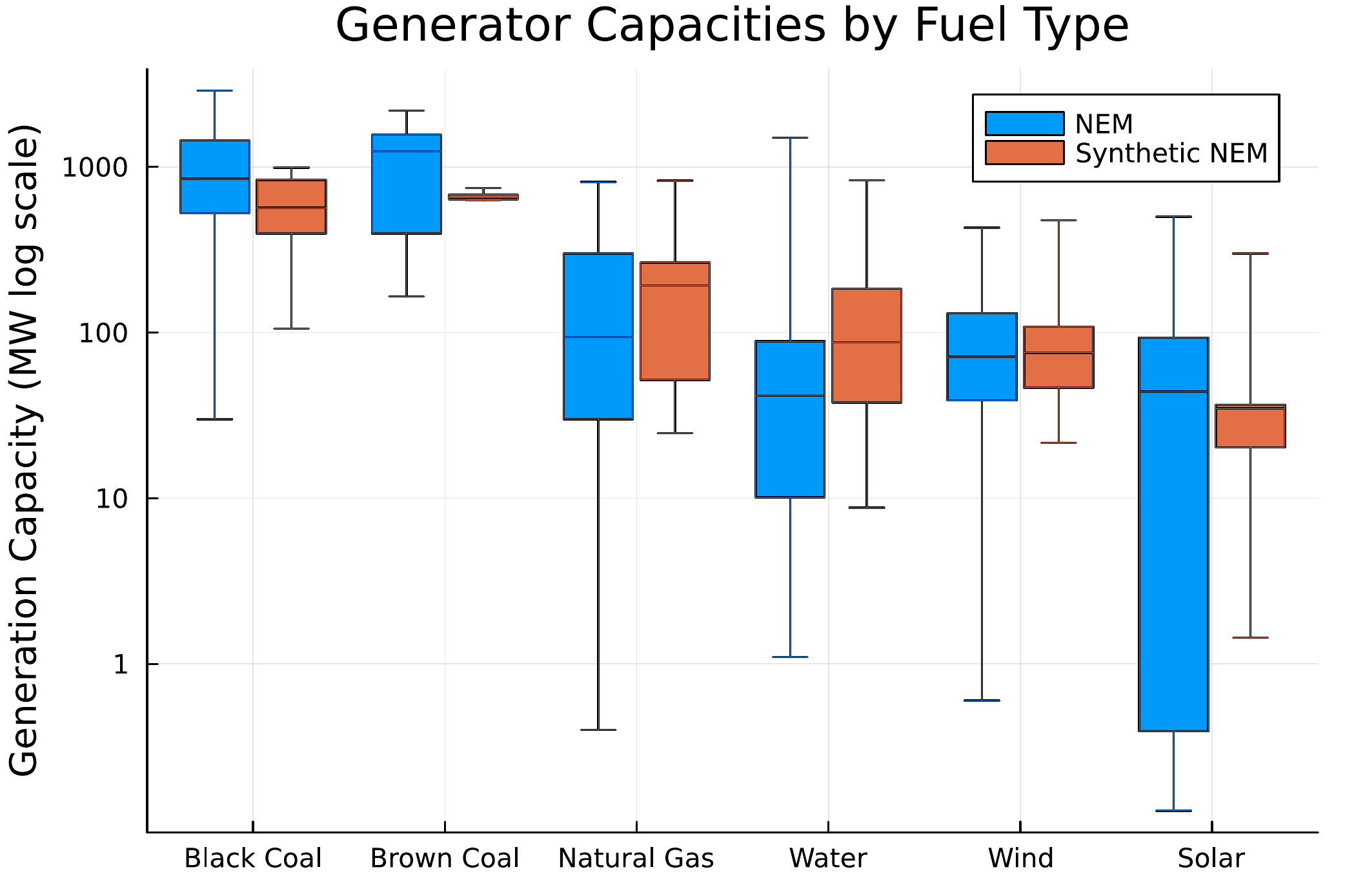}
  \caption{Comparing NEM and S-NEM2000 w.r.t generator count by fuel type.}
  \label{fig:gen_compare_fueltypes_by_capacity}
\end{figure}

Figure~\ref{fig:gen_compare_fueltypes_by_capacity} shows comparison of the generator capacities categorized by fuel type between the NEM and the S-NEM2000.  Quantiles comparison shows that the classification has rendered a reasonable traction of generator capacities for each fuel type.

\begin{figure}[H]
  \centering
  \includegraphics[width=0.9\columnwidth]{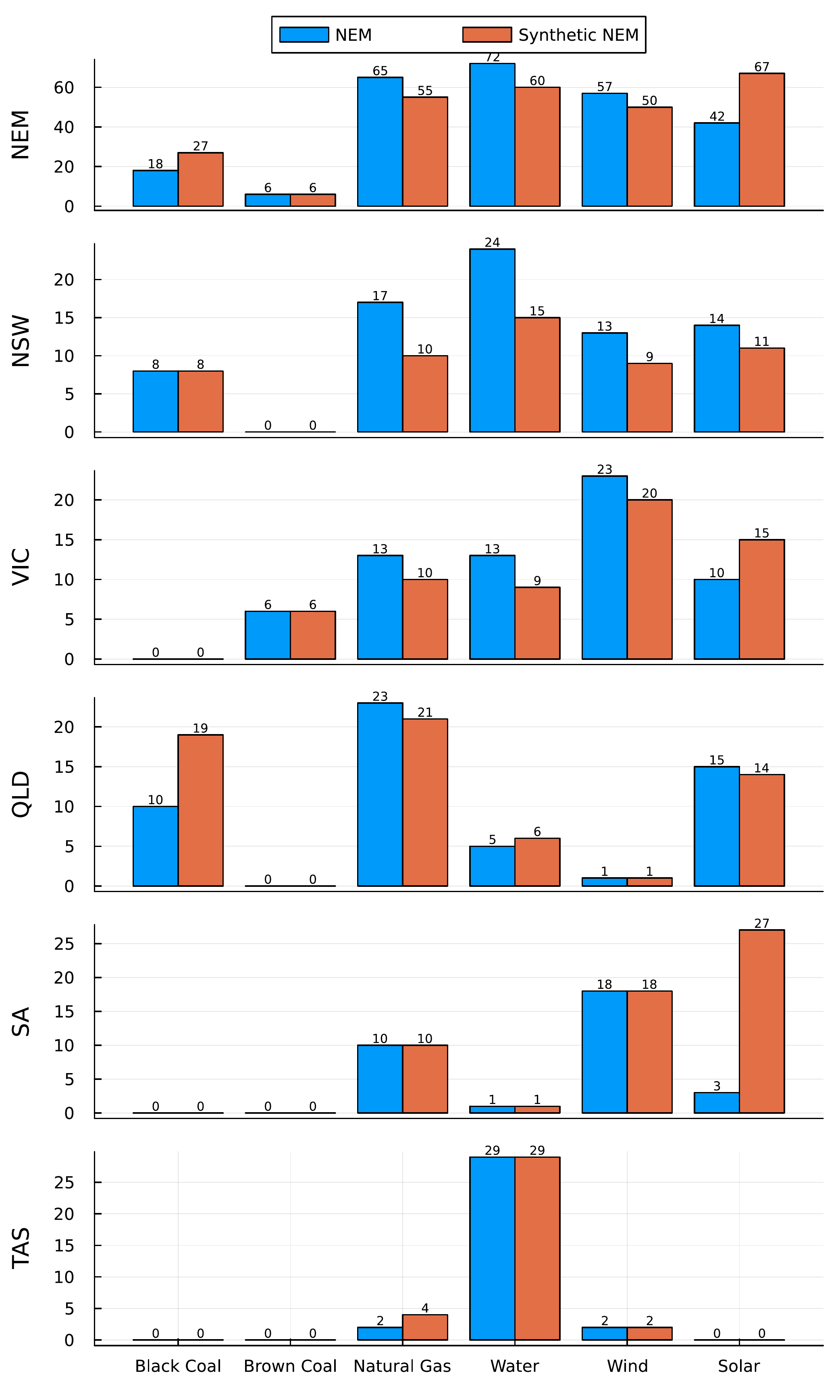}
  \caption{Comparing NEM and S-NEM2000 w.r.t generator count by fuel type and state.}
  \label{fig:gen_compare_fueltypes_by_count_by_state}
\end{figure}

\begin{figure}[H]
  \centering
  \includegraphics[width=0.9\columnwidth]{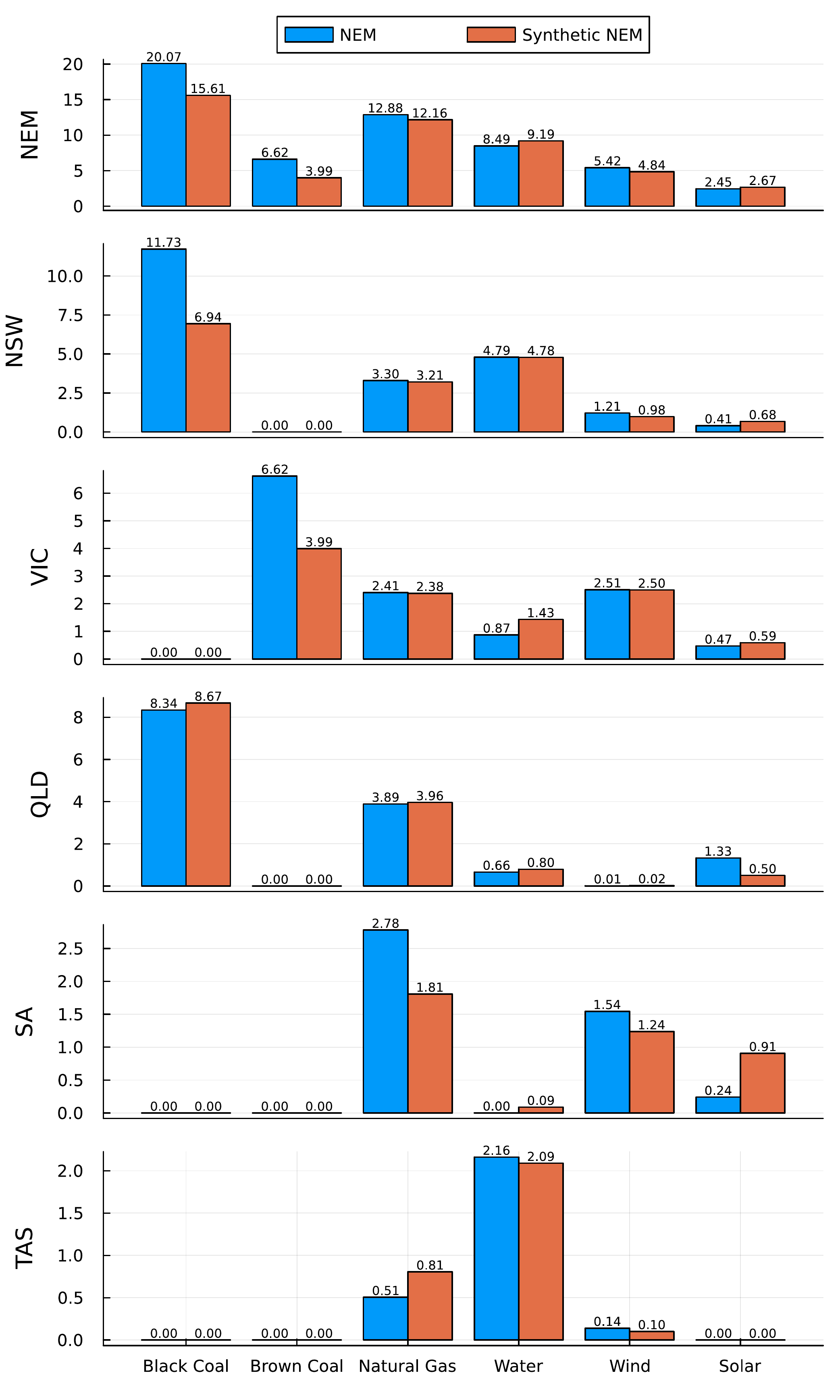}
  \caption{Comparing NEM and S-NEM2000 w.r.t generator capacity by fuel type and state.}
  \label{fig:gen_compare_fueltypes_by_capacity_by_state}
\end{figure}

To further assess the performance of classification, the predicted generation type was used to aggregate the generation capacity and the count of the synthetic generators in each state to compare to the aggregate generation capacity and count of generators within the NEM \cite{GA2021}.  These aggregate comparisons can be seen in Figures \ref{fig:gen_compare_fueltypes_by_count_by_state} and \ref{fig:gen_compare_fueltypes_by_capacity_by_state}. 
These figures show how the total generation capacity and count of each state compare between the classified synthetic model and the reference data.  The aggregate generation capacity and count categories in both network models follow a similar pattern while the synthetic network has a lower total capacity level. The shortage generation is because the reference data is an updated list of all existing generators in the NEM, while the S-NEM2300 model was developed from a PSS/E model of the network representing a 2018 Summer day. The synthetic model was developed to replicate that specific moment in time, which caused the exclusion of several generation units that were not in operation at the time. It is therefore realistic to see a lower generation capacity level in the synthetic network model. 
More observations from the Figures~\ref{fig:gen_compare_fueltypes_by_count_by_state} and \ref{fig:gen_compare_fueltypes_by_capacity_by_state} are as follows:
\begin{itemize}
	\item The total generation capacity and generators count for the whole NEM, classified by fuel type,  are often lower in the S-NEM2000,  except for solar which is due to consideration of all `other' generation units such as var compensators and battery storage.
	\item Tasmania generation capacity and count are represented quite well. 
	\item Queensland is represented reasonably well in terms of generator count and capacity.
	\item Generation capacity and count in New South Wales, Victoria and South Australia have been represented fairly reasonable for most fuel types and less accurate for some others.
\end{itemize}

Overall, the classification demonstrates a reasonable representation of fuel types in each state and across the whole NEM. One could also further modify generator capacities to match the actual aggregated value to be a closer representation of the actual system.

\subsection{Generation Cost and Operation Data}
After generator fuel type classification we use AEMO's Fuel and Technology Cost Review \cite{AEMO2014} to assign cost and operational models for different generators. Table~\ref{table_generators_models} reports the main parameters identified from the reference file and assigned to the fuel categories.

\begin{table}[tbh]
  \centering 
   \caption{S-NEM2000 generator cost and operational data.} 
   \label{table_generators_models}
    \begin{tabular}{l | l l l l l l}
    \hline
       & Black coal & Brown coal & Natural Gas & Water &  Wind & Solar \\
    \hline
Min Gen ($\%$ of nameplate capacity) & 40 & 60 & 0 & 0 & 0 & 10 \\
Minimum On/Off Time (Hours)	& 8/8 & 16/16 & 1/1 & 1/1 & 1/1 & 1/1 \\
No Load Fuel Consumption \\ \quad ($\%$ of Full Load Fuel Consumption) & 10 & 10 & 30 & 0 & 0 & 0 \\
Auxiliary load \\ \quad ($\%$ of nameplate capacity) & 6 & 10 & 3 & 1 & 1 & 8 \\
Ramp Up Rate (MW/h) & 120 & 130 & 600 & 7630 & 900 & 10000 \\
Ramp Down Rate (MW/h) & 150 & 140 & 600 & 1840 & 600 & 10000 \\
Thermal Efficiency ($\%$, HHV sent-out)	& 35.4 & 23.5 & 32 & 100 & 100 & 100 \\
Maintenance (days/yr) & 20 & 20 & 5 & 7 & 3 & 10 \\
Fixed Op Cost ($\$$/MW/year) & 53400 & 147200 & 14200 & 56700 & 43000 & 64000 \\
Variable Op Cost ($\$$/MWh sent-out) & 1 & 1 & 8.3 & 7 & 10 & 15.2 \\
No Load Costs ($\$$/MW, 2014-15) & 1500 & 1500 & 500 & 500 & 2000 & 500	 \\
No Load Recurring Cost ($\$$/MW/year) & 500 & 500 & 200 & 200 & 1000 & 200 \\
Cold/Warm/Hot Start-up Costs ($\$$/MW) & 350 & 350 & 100 & 5 & 5 & 5 \\
Warm Start-up Costs ($\$$/MW) & 120 & 120 & 100 & 3 & 3 & 3	 \\
Hot Start-up Costs ($\$$/MW) & 40 & 40 & 100 & 2 & 2 & 2	 \\

\hline
\end{tabular}
\end{table}

\section{Proof of Concept Study}
\label{sec:proof_of_concept}
In addition to the integrated NEM data (\texttt{snem2000.m}), we divide it into two sub-networks 1) mainland (\texttt{snem1803.m}) and 2) Tasmania (\texttt{snem197.m}).
In this section we run optimal power flow on S-NEM2000 and its mainland and Tasmanis sub-networks for AC power flow relaxations and compare objective values and runtime, which is similar to how the results in \cite{PGlib2021} are presented. We use PowerModels.jl version 0.19.8 and the interior point solver Ipopt \cite{Wachter2006} with MUMPS linear solver to solve the optimal power flow problems. The computation is performed on a PC with  6-Core Intel(R) Core(TM) i7 @ 2.60GHz processor with 32\,GB of RAM.

\begin{table}[htpb]
    \centering
     \caption{Quality and runtime results of AC power flow formulations.}
    \begin{tabular}{c c c c c}
    \\
    Network & Formulation & Objective value  ($\$$) &  Gap ($\%$) & Runtime (s) \\
    \hline
    Whole NEM& NLP  & 98322.1 & - & 19.75 \\
    (\texttt{snem2000.m}) & SOC  & 90424.6 & 8.03 & 18.83 \\
    & DC   & 87396.6 & 11.11 & 1.55 \\
    \hline
    Mainland & NLP  & 98320.6 & - & 17.44 \\
    (\texttt{snem1803.m}) & SOC  & 90423.1 & 8.03 & 16.54 \\
    & DC   & 87395.0 & 11.11 & 2.54 \\
    \hline
    Tasmania& NLP  & 1.50165 & - & 0.64 \\
    (\texttt{snem197.m}) & SOC  & 1.50067 & 0.06 & 0.72 \\
    & DC   & 1.47406 & 1.84 & 0.06 \\
    \hline
    \end{tabular}
    \label{tab:opf_formulations}
\end{table}

\section{Conclusion and Future Work}
In this study, we utilized an existing synthetic representation of the NEM and developed the S-NEM2000 benchmark for power system optimization studies. By converting the synthetic network data into the PowerModels.jl and MATPOWER data models, we enabled feasibility and validation checks through power flow and optimal power flow studies. The resulting benchmark model includes enhanced features such as thermal limits and generation fuel type and cost models, providing a valuable resource for optimization studies in the power system domain.
We release the whole of NEM data (\texttt{snem2000.m}) under the creative commons `CC BY 4.0'\footnote{\url{https://creativecommons.org/licenses/by/4.0/}} license at \url{https://github.com/csiro-energy-systems/Synthetic-NEM-2000bus-Data}. The authors also intend to publish this data set in Power Grid Library Optimal Power Flow benchmarks.

However, it is important to note that the synthetic model used in this study is not an exact replica of the actual NEM, and there is still room for improvement. Future work should focus on incorporating additional components, such as High-Voltage Direct Current (HVDC) lines, which play a significant role in interconnecting regions within the NEM. The inclusion of HVDC lines would enable more accurate simulations and analyses of power flow and optimization scenarios.
Furthermore, expanding the benchmark model to include detailed load and generation profiles is another important avenue for future work. Incorporating realistic and time-varying data on load and generation patterns would provide a more comprehensive representation of the NEM and enable more accurate assessments of system performance and optimization strategies. Additionally, for the data-driven generation classification approach, using additional generator properties, for example voltage, could improve the classification results. 

The ongoing development and growth of this work will involve continuous refinement and expansion of the S-NEM2000 model, incorporating additional components and data sources to enhance its representation of the real-world system. These advancements will contribute to furthering our understanding of the NEM and supporting effective decision-making in the context of power system optimization and planning.

\section*{Acknowledgement}
We would like to express our gratitude to Dr. Felipe Arraño-Vargas and Dr. Georgios Konstantinou who released the original data used in this study, and for their ongoing support and advice on gaining insights from the data. Their invaluable contribution and commitment to open data have facilitated our research and enabled us to develop the S-NEM2000 benchmark for power system optimization studies.

We extend our appreciation to AEMO and GeoScience Australia for providing the publicly available generation datasets that formed the basis of our analysis. Their efforts in collecting and maintaining comprehensive data on the NEM have been instrumental in enhancing our understanding of the power system domain and played a crucial role in assigning appropriate properties to the synthetic generators in our model and ensuring its fidelity to the NEM.

\bibliographystyle{IEEEtran}
\bibliography{library.bib}

\end{document}